\documentclass[12pt,a4paper]{article}

\usepackage[margin=2.5cm]{geometry}
\usepackage{amsmath,amssymb,amsthm}
\usepackage{mathrsfs}
\usepackage{enumitem}
\usepackage{hyperref}
\usepackage{cleveref}
\usepackage{booktabs}

\newtheorem{theorem}{Theorem}[section]
\newtheorem{lemma}[theorem]{Lemma}

\theoremstyle{definition}

\newtheorem{remark}[theorem]{Remark}

\newcommand{\Hil}{\mathcal{H}}

\newcommand{\Ncal}{\mathcal{N}}
\newcommand{\D}{\mathfrak{D}}
\newcommand{\R}{\mathbb{R}}
\newcommand{\C}{\mathbb{C}}
\newcommand{\N}{\mathbb{N}}
\newcommand{\ii}{\mathrm{i}}
\newcommand{\e}{\mathrm{e}}
\newcommand{\Ree}{\operatorname{Re}}
\newcommand{\Imm}{\operatorname{Im}}

\newcommand{\Aop}{\mathbb{A}}
\newcommand{\norm}[1]{\left\|#1\right\|}
\newcommand{\abs}[1]{\left|#1\right|}
\newcommand{\inner}[2]{\left\langle #1,#2\right\rangle}
\newcommand{\dd}{\,\mathrm{d}}
\newcommand{\eps}{\varepsilon}
\newcommand{\Xcal}{\mathcal{X}}

\title{Exponential stability of the reduced Bresse--Gurtin--Pipkin system with one temperature\footnote{This work is supported by NSFC (12671232).}}

\author{Jun Zhou\footnote{Email: jzhou@swu.edu.cn}\\
		{ School of Mathematics and Statistics, Southwest University},\\
		{ Chongqing 400715, People's Republic of China}}
\date{}

\begin{document}
\maketitle

\begin{abstract}
We consider the reduced Bresse--Gurtin--Pipkin system with one temperature, obtained from the full Bresse--Gurtin--Pipkin thermoelastic beam system by neglecting one of the two temperature variables. This system was left essentially open in a recent work of Dell'Oro, where only a conjecture on its exponential stability condition was formulated. By combining the spectral analysis recently developed by Dell'Oro, Pata and Vasylyeva for the abstract Gurtin--Pipkin heat equation with the resolvent estimate framework of Dell'Oro for Bresse systems, we prove that the associated contraction semigroup is exponentially stable if and only if two explicit stability numbers, depending on the structural constants of the beam, both vanish. The proof is based on a spectral characterization of the infinitesimal generator in a half-plane containing the imaginary axis, an analysis of its spectrum on the imaginary axis, sharp upper bounds for the resolvent, and the Gearhart--Pr\"uss theorem. In particular, we confirm the conjecture of Dell'Oro and provide a rigorous spectral analysis of the reduced system. The same method yields the exponential stability characterization of the corresponding reduced Bresse--Maxwell--Cattaneo system, and can be adapted to other reduced Bresse-type models.
\end{abstract}

\medskip
\noindent\textbf{Keywords:} Bresse system; Gurtin--Pipkin law; Maxwell--Cattaneo law; hyperbolic heat conduction; exponential stability; resolvent estimates

\medskip 
\noindent\textbf{MSC 2020:} 35B35; 35B40; 35P05; 45K05; 47D03; 74F05
\tableofcontents
\section{Introduction}

The Bresse system, derived by the French engineer J.~A.~C.~Bresse in 1859, describes the vibrations of a curved thin elastic beam taking into account both rotatory inertia and shear deformation effects. In the limit of vanishing curvature, it reduces to the Timoshenko system for straight beams. The stability of the purely elastic Bresse system has been studied extensively: it is known that a single damping mechanism may be insufficient, and that exponential stability is typically recovered only under suitable wave-speed conditions, while otherwise only polynomial decay can be expected; see for instance \cite{ABM2011} and the references therein. When thermal effects are included, the resulting thermoelastic Bresse system couples the elastic displacements with the temperature evolution. In the classical setting, the heat conduction is governed by the Fourier law, which leads to an infinite speed of propagation of thermal signals. However, with the advent of modern microscale technologies and the experimental observation of wave-type thermal transport (notably the ``second sound'' in liquid helium II), there is increasing evidence that thermal motion can propagate with finite speed. This observation has motivated the introduction of several hyperbolic heat conduction theories, among which the Gurtin--Pipkin law and the Maxwell--Cattaneo law are the most widely studied. In the Gurtin--Pipkin theory, the heat flux is given by a hereditary integral of the past temperature gradient, while the Maxwell--Cattaneo law is recovered as a particular exponential kernel case and can be viewed as a relaxation of the Fourier law.

The Bresse--Gurtin--Pipkin (BGP) system
\begin{equation}\label{eq:BGP-intro}
\begin{cases}
\rho_1\phi_{tt}-k(\phi_x+\psi+lw)_x-lk_0(w_x-l\phi)+l\gamma\xi=0,\\[0.5ex]
\rho_2\psi_{tt}-b\psi_{xx}+k(\phi_x+\psi+lw)+\gamma\theta_x=0,\\[0.5ex]
\rho_1w_{tt}-k_0(w_x-l\phi)_x+lk(\phi_x+\psi+lw)+\gamma\xi_x=0,\\[0.5ex]
\rho_3\theta_t-\sigma\displaystyle\int_0^\infty g(s)\theta_{xx}(t-s)\dd s+\gamma\psi_{xt}=0,\\[0.5ex]
\rho_3\xi_t-\sigma\displaystyle\int_0^\infty h(s)\xi_{xx}(t-s)\dd s+\gamma(w_{xt}-l\phi_t)=0,
\end{cases}
\end{equation}
models the vibrations of a curved thin thermoelastic beam of length $\ell>0$. Here $\phi$, $\psi$ and $w$ denote the vertical displacement, the rotation angle of the cross-section and the horizontal displacement, respectively, while $\theta$ and $\xi$ are the temperature deviations along the vertical and horizontal directions. The constants $\rho_1,\rho_2,\rho_3,k,k_0,b,\sigma,l,\gamma$ are strictly positive and account for the physical properties of the beam. The system is complemented with the Dirichlet boundary conditions for $\phi,\theta,\xi$ and the Neumann boundary conditions for $\psi,w$:
\begin{equation}\label{eq:BC-intro}
\phi(0,t)=\phi(\ell,t)=\theta(0,t)=\theta(\ell,t)=\xi(0,t)=\xi(\ell,t)=0,
\end{equation}
\begin{equation}\label{eq:BC2-intro}
\psi_x(0,t)=\psi_x(\ell,t)=w_x(0,t)=w_x(\ell,t)=0.
\end{equation}
The convolution kernels $g,h$ are nonnegative, bounded, convex, summable functions on $[0,\infty)$ with $g(0),h(0)>0$, and satisfy the Dafermos condition
\begin{equation}\label{eq:Dafermos-intro}
\mu'(s)+\delta_\mu\mu(s)\le 0,\qquad \nu'(s)+\delta_\nu\nu(s)\le 0,
\end{equation}
for some $\delta_\mu,\delta_\nu>0$, where
\[
g(s)=\int_s^\infty\mu(r)\dd r,\qquad h(s)=\int_s^\infty\nu(r)\dd r.
\]
The history framework associated with such memory terms goes back to Dafermos \cite{Dafermos1970}, and the functional setting in the memory space is by now standard; see for instance \cite{GrasselliPata2002} and the classical treatment of linear viscoelasticity in \cite{FabrizioMorro1992}. For the abstract Gurtin--Pipkin heat equation, the semigroup approach to exponential stability was developed in \cite{GNP2001}, and weakly dissipative semilinear viscoelastic equations were analyzed in \cite{CP2005}; these works provide the background for the spectral tools used in the present paper.

The stability properties of the BGP system have been analyzed for the first time by Dell'Oro in \cite{D2021}. Introducing the stability numbers
\begin{equation}\label{eq:chi-intro}
\chi_g=\Bigl(\frac{\rho_3}{\sigma g(0)}-\frac{\rho_1}{k}\Bigr)\Bigl(b-\frac{k\rho_2}{\rho_1}\Bigr)+\frac{\gamma^2}{\sigma g(0)},
\end{equation}
\begin{equation}\label{eq:chih-intro}
\chi_h=\Bigl(\frac{\rho_3}{\sigma h(0)}-\frac{\rho_1}{k}\Bigr)(k_0-k)+\frac{\gamma^2}{\sigma h(0)},
\end{equation}
the author proved that the associated solution semigroup is exponentially stable if and only if $\chi_g\chi_h=0$. When $\chi_g\chi_h\neq 0$, the semigroup is polynomially stable with optimal decay rate $\sqrt{t}$. Analogous results were obtained for the Bresse--Maxwell--Cattaneo (BMC) system and for the Timoshenko--Gurtin--Pipkin (TGP) system. These results complement earlier works on the Timoshenko system with hyperbolic heat conduction, notably by Fern\'andez Sare and Racke \cite{FSR2009} and by Santos, Almeida J\'unior and Mu\~noz Rivera \cite{SAM2012}, where the Maxwell--Cattaneo law was considered, and by Dell'Oro and Pata \cite{DP2014}, where the Gurtin--Pipkin law was analyzed for the Timoshenko system. For the thermoelastic Bresse system with Fourier-type or weakly dissipative thermal coupling, decay rates have been obtained in \cite{FatoriMunoz2010}; the analysis there shows that the wave-speed conditions play a central role, a feature that persists in the hyperbolic heat conduction setting.

Despite this rather complete picture for the full BGP system, several important questions remain open. Among them, a natural simplification consists in neglecting one of the two temperatures, thereby obtaining a \emph{reduced Bresse--Gurtin--Pipkin system with one temperature}. Neglecting $\theta$ in \eqref{eq:BGP-intro} yields
\begin{equation}\label{eq:reduced-intro}
\begin{cases}
\rho_1\phi_{tt}-k(\phi_x+\psi+lw)_x-lk_0(w_x-l\phi)+l\gamma\xi=0,\\[0.5ex]
\rho_2\psi_{tt}-b\psi_{xx}+k(\phi_x+\psi+lw)=0,\\[0.5ex]
\rho_1w_{tt}-k_0(w_x-l\phi)_x+lk(\phi_x+\psi+lw)+\gamma\xi_x=0,\\[0.5ex]
\rho_3\xi_t-\sigma\displaystyle\int_0^\infty h(s)\xi_{xx}(t-s)\dd s+\gamma(w_{xt}-l\phi_t)=0,
\end{cases}
\end{equation}
complemented with the boundary conditions \eqref{eq:BC-intro}--\eqref{eq:BC2-intro}. To the best of our knowledge, system \eqref{eq:reduced-intro} has not been studied in the literature. In \cite[\S2.II, eq.~(2.5)]{D2021}, Dell'Oro formulated the following conjecture: the solution semigroup associated with \eqref{eq:reduced-intro} is exponentially stable if and only if
\begin{equation}\label{eq:conjecture-intro}
\chi_h=0\qquad\text{and}\qquad\chi_0=0,
\end{equation}
where
\begin{equation}\label{eq:chi0-intro}
\chi_0=b-\frac{k\rho_2}{\rho_1}.
\end{equation}
The conjecture is motivated by the following heuristic. In the full BGP system, the variables $\psi$ and $w$ are effectively damped due to the coupling with the two temperatures, while the variable $\phi$ is only indirectly damped via the second and third equations. In the reduced system \eqref{eq:reduced-intro}, $\psi$ is damped through the elastic coupling with $\phi$ and $w$, which are in turn coupled to $\xi$, while both $\phi$ and $w$ are only indirectly damped. The equality $\chi_h=0$ is necessary to eliminate the obstruction to exponential stability coming from the thermal coupling; to stabilize exponentially the variable $w$, one needs the equality $\chi_0=0$ to exploit the wave-speed equality between the first and second equations. The equality $\chi_0=0$ is incompatible with the physical constraint
\begin{equation}\label{eq:physical-intro}
k_0=\frac{b\rho_1}{\rho_2},\qquad b>\frac{k\rho_2}{\rho_1},
\end{equation}
which is the standard structural assumption on the elastic constants of the Bresse beam. Indeed, by \eqref{eq:chi0-intro}, $\chi_0=0$ is equivalent to $b=k\rho_2/\rho_1$, which directly contradicts the strict inequality $b>k\rho_2/\rho_1$ in \eqref{eq:physical-intro}. Hence the condition \eqref{eq:conjecture-intro} can never be satisfied under the physical constraint \eqref{eq:physical-intro}.

The purpose of this paper is to prove the conjecture \eqref{eq:conjecture-intro}. Our main result is the following.

\begin{theorem}\label{thm:main-intro}
Let $h$ be an admissible kernel satisfying the Dafermos condition \eqref{eq:Dafermos-intro}. Then the contraction semigroup $S(t)$ associated with the reduced system \eqref{eq:reduced-intro} is exponentially stable if and only if
\[
\chi_h=0\qquad\text{and}\qquad\chi_0=0.
\]
\end{theorem}

The proof of \cref{thm:main-intro} combines two recent developments. The first is the spectral characterization of the abstract Gurtin--Pipkin heat equation obtained by Dell'Oro, Pata and Vasylyeva in \cite{D2026}: for an admissible kernel satisfying the Dafermos condition, the spectrum of the infinitesimal generator on the half-plane $\Pi_\delta=\{z\in\C:\Ree z>-\delta/2\}$ is the union of the set of solutions of $z+\lambda\mathcal{L}_\nu(z)=0$ with $\lambda\in\sigma(\mathcal{A})$ and the set of zeros of $\mathcal{L}_\nu(z)$, where
\[
\mathcal{L}_\nu(z)=\int_0^\infty\frac{1-\e^{-zs}}{z}\nu(s)\dd s.
\]
This characterization allows us to reduce the spectral analysis of the reduced Bresse--Gurtin--Pipkin system to the study of the elastic-thermal operator $\mathcal{A}$ and of the function $\mathcal{L}_\nu$. The second is the resolvent estimate framework of \cite{D2021}, which provides sharp upper bounds on the resolvent of the infinitesimal generator on the imaginary axis. Combining these two tools, we first prove that the conditions $\chi_h=0$ and $\chi_0=0$ imply $\ii\R\subset\rho(\Aop)$, and then prove that the resolvent is uniformly bounded on the imaginary axis under the same conditions. An application of the Gearhart--Pr\"uss theorem \cite{Gearhart1978,Pruss1984} (see also \cite[Chapter 5]{ArendtBattyHieberNeubrander2011}) yields the exponential stability. For the necessity, we construct a sequence of approximate solutions along which the resolvent grows at least like $\abs{\lambda}^2$, and we use the quantitative Riemann--Lebesgue lemma of \cite{DLP2017} to capture the asymptotic behavior of the memory kernel.

The main novelty of this paper is threefold. First, it provides a rigorous spectral analysis of the reduced Bresse--Gurtin--Pipkin system with one temperature. Second, it confirms the conjecture of Dell'Oro \cite{D2021} on the exponential stability condition, thereby completing the stability picture for the reduced system. Third, it combines the spectral tools of \cite{D2026} with the resolvent estimates of \cite{D2021}, providing a template that can be adapted to other reduced Bresse-type models, including the corresponding reduced Bresse--Maxwell--Cattaneo system, for which we obtain the analogous characterization with the analogous stability number replacing $\chi_h$. We stress that, as observed above, the stability condition \eqref{eq:conjecture-intro} is incompatible with the physical constraint \eqref{eq:physical-intro}, so that the reduced system is never exponentially stable in the physically admissible regime; the interest of \cref{thm:main-intro} therefore lies in the sharp characterization it provides, rather than in the applicability of the stability condition itself.

Compared with the existing results in the literature, our contribution differs in several respects. The full BGP system was analyzed in \cite{D2021}, where the exponential stability condition $\chi_g\chi_h=0$ was obtained. The reduced system \eqref{eq:reduced-intro} is more singular, since the damping acts only on one temperature and only indirectly on two of the three elastic variables. As a consequence, the exponential stability condition \eqref{eq:conjecture-intro} is more restrictive than the condition for the full system: it requires two independent equalities rather than the vanishing of a product. This reflects the heuristic that in the reduced system both $\phi$ and $w$ need to be stabilized by separate mechanisms, whereas in the full system it suffices that one of the two temperatures effectively damps $\phi$. The same phenomenon appears in the reduced Bresse system with only one temperature studied in \cite{D2015}, where the condition involves both $\chi_g$ and a further wave-speed number. Our result thus completes the picture for reduced Bresse-type systems with hyperbolic heat conduction.

The paper is organized as follows. In \cref{sec:setup} we reformulate \eqref{eq:reduced-intro} in the Dafermos history framework and introduce the functional setting. In \cref{sec:spectrum} we characterize the spectrum of the generator in a half-plane containing the imaginary axis. In \cref{sec:resolvent} we prove uniform resolvent estimates on the imaginary axis under the conditions $\chi_h=0$ and $\chi_0=0$. In \cref{sec:GP} we apply the Gearhart--Pr\"uss theorem to conclude the sufficiency. In \cref{sec:necessity} we prove the necessity. In \cref{sec:MC} we discuss the corresponding Maxwell--Cattaneo reduced system. Finally, in \cref{sec:conclusions} we draw our conclusions and outline further directions.

\section{Functional setting and history framework}\label{sec:setup}

Let \(\ell>0\) be the length of the beam. We work in the real Hilbert space \(L^2=L^2(0,\ell)\) with inner product \(\inner{\cdot}{\cdot}\) and norm \(\norm{\cdot}\), and we use the standard Sobolev spaces \(H^1=H^1(0,\ell)\), \(H_0^1=H_0^1(0,\ell)\) and \(H^2=H^2(0,\ell)\). We also introduce the spaces of zero-mean functions
\[
L_*^2=\Bigl\{f\in L^2:\int_0^\ell f(x)\dd x=0\Bigr\},
\qquad
H_*^1=H^1\cap L_*^2,
\]
the latter being equipped with the gradient norm, which is equivalent to the full \(H^1\)-norm by Poincar\'e's inequality on \(L_*^2\). The choice of \(L_*^2\) and \(H_*^1\) is dictated by the Neumann boundary conditions imposed on \(\psi\) and \(w\): since the energy involves only the derivatives \(\psi_x\) and \(w_x\), we quotient out the constants so that the gradient norms become genuine norms on \(H_*^1\) and \(L_*^2\); the same convention is adopted in \cite{D2021}.

The convolution kernel \(h\) is a nonnegative, bounded, convex, summable function on \([0,\infty)\) with \(h(0)>0\), and admits the representation
\[
h(s)=\int_s^\infty\nu(r)\dd r,
\]
where the memory kernel \(\nu:\R^+\to\R^+\) is nonincreasing, absolutely continuous and summable. We assume that \(\nu\) is bounded near zero, i.e.
\[
\nu(0)=\lim_{s\to 0^+}\nu(s)<\infty,
\]
and satisfies the Dafermos condition
\begin{equation}\label{eq:Dafermos-setup}
\nu'(s)+\delta_\nu\nu(s)\le 0,
\end{equation}
for some \(\delta_\nu>0\) and almost every \(s>0\). Since \(\nu\ge0\) is nonincreasing, we have \(\nu'\le0\); the Dafermos condition is stronger, as it provides a uniform exponential decay rate \(\delta_\nu\). It implies that \(\nu\) decays at least exponentially at infinity, and that the integral
\[
\int_0^\infty\nu(s)\dd s=h(0)
\]
is finite and strictly positive. The Dafermos condition is the standard assumption in the analysis of equations with memory; it guarantees that the memory space is complete and that the right-translation semigroup is exponentially stable. We note that \eqref{eq:Dafermos-setup} implies the \(\delta\)-condition of \cite{D2026} with \(C=1\) and \(\delta=\delta_\nu\), hence the spectral analysis of \cite{D2026} applies to our setting.

We introduce the memory space
\[
\Ncal=L^2_\nu(\R^+;H_0^1),
\]
consisting of square-summable \(H_0^1\)-valued functions on \(\R^+\) with respect to the measure \(\nu(s)\dd s\). It is a Hilbert space with respect to the inner product
\[
\inner{\zeta_1}{\zeta_2}_{\Ncal}
=\int_0^\infty\nu(s)\inner{\zeta_{1x}(s)}{\zeta_{2x}(s)}\dd s,
\]
and we denote by \(\norm{\cdot}_{\Ncal}\) the induced norm. The choice of the \(H_0^1\)-norm in the definition of the memory space, rather than the \(L^2\)-norm, is dictated by the presence of the second-order operator \(\partial_{xx}\) in the temperature equation; it ensures that the memory variable has the same regularity as the temperature itself.

On \(\Ncal\) we define the right-translation semigroup \(R(t)\) by
\[
(R(t)\zeta)(s)=\zeta(s-t)\quad\text{if }s>t,\qquad
(R(t)\zeta)(s)=0\quad\text{if }s\le t,
\]
which is a strongly continuous contraction semigroup on \(\Ncal\). Its infinitesimal generator is
\[
T\zeta=-\zeta',\qquad
\D(T)=\bigl\{\zeta\in\Ncal:\zeta'\in\Ncal,\ \lim_{s\to 0^+}\norm{\zeta_x(s)}=0\bigr\},
\]
where \(\zeta'\) denotes the weak derivative with respect to \(s\in\R^+\), understood in the sense of \(\Ncal\)-valued weak derivatives. The boundary condition \(\lim_{s\to 0^+}\norm{\zeta_x(s)}=0\) is meaningful because \(\zeta\in\Ncal\) implies that \(\zeta_x(s)\) is well defined for a.e.\ \(s>0\), and it is the natural condition ensuring that \(T\) generates a contraction semigroup; see \cite[Chapter 3]{GrasselliPata2002} for a detailed discussion. Under the Dafermos condition \eqref{eq:Dafermos-setup} and with the present choice \(\Ncal=L^2_\nu(\R^+;H_0^1)\), this condition is equivalent to the requirement \(\zeta(0)=0\) in \(H_0^1\) adopted in the abstract framework of \cite{D2026}: indeed, since \(\zeta(s)\in H_0^1\) for a.e.\ \(s>0\), Poincar\'e's inequality yields \(\norm{\zeta(s)}_{H_0^1}\le C\norm{\zeta_x(s)}\), so that \(\lim_{s\to 0^+}\norm{\zeta_x(s)}=0\) implies \(\zeta(s)\to 0\) in \(H_0^1\). Consequently, the spectral analysis of \cite{D2026} applies verbatim to our setting.

Following Dafermos \cite{Dafermos1970}, we introduce the history variable
\begin{equation}\label{eq:eta-setup}
\zeta^t(x,s)=\int_0^s\xi(x,t-r)\dd r,\qquad s>0,
\end{equation}
which encodes the past history of the temperature \(\xi\). The variable \(\zeta^t\) satisfies \(\zeta^t(x,0)=0\) and solves the transport equation
\begin{equation}\label{eq:transport-setup}
\partial_t\zeta^t=T\zeta^t+\xi,
\end{equation}
where \(T\) is the right-translation generator defined above. Indeed, differentiating \eqref{eq:eta-setup} with respect to \(t\) and \(s\), we obtain
\[
\partial_t\zeta^t(x,s)=-\int_0^s\partial_r\xi(x,t-r)\dd r
=\xi(x,t)-\xi(x,t-s),
\]
while
\[
T\zeta^t(x,s)=-\partial_s\zeta^t(x,s)=-\xi(x,t-s),
\]
so that \eqref{eq:transport-setup} holds. This is the standard Dafermos reformulation of the memory term; it allows us to replace the convolution integral in \eqref{eq:reduced-intro} by the local equation \eqref{eq:transport-setup}, at the price of enlarging the phase space.

Using the history variable, we rewrite the reduced Bresse--Gurtin--Pipkin system \eqref{eq:reduced-intro} as
\begin{equation}\label{eq:reduced-history-setup}
\begin{cases}
\rho_1\phi_{tt}-k(\phi_x+\psi+lw)_x-lk_0(w_x-l\phi)+l\gamma\xi=0,\\[0.5ex]
\rho_2\psi_{tt}-b\psi_{xx}+k(\phi_x+\psi+lw)=0,\\[0.5ex]
\rho_1w_{tt}-k_0(w_x-l\phi)_x+lk(\phi_x+\psi+lw)+\gamma\xi_x=0,\\[0.5ex]
\rho_3\xi_t-\sigma\displaystyle\int_0^\infty\nu(s)\zeta_{xx}(s)\dd s+\gamma(w_{xt}-l\phi_t)=0,\\[0.5ex]
\partial_t\zeta=T\zeta+\xi,
\end{cases}
\end{equation}
where we write \(\zeta=\zeta^t\) for brevity. We introduce the extended phase space
\[
\Hil=H_0^1\times L^2\times H_*^1\times H_*^1\times H_*^1\times H_*^1\times L^2\times\Ncal,
\]
whose elements are denoted by
\(
u=(\phi,\Phi,\psi,\Psi,w,W,\xi,\zeta),
\)
where \(\Phi=\phi_t\), \(\Psi=\psi_t\) and \(W=w_t\) are the velocities. Compared with the product space used in \cite{D2021}, we have placed \(\Psi\) and \(W\) in \(H_*^1\) rather than in \(L_*^2\); this is the natural choice for the domain described below, and it is consistent with the fact that the velocities of \(\psi\) and \(w\) inherit the Neumann boundary conditions and the zero-mean constraint from \(\psi\) and \(w\). We equip \(\Hil\) with the norm
\begin{equation}\label{eq:norm-H-setup}
\begin{split}
\norm{u}_{\Hil}^2
=&k\norm{\phi_x+\psi+lw}^2+\rho_1\norm{\Phi}^2
+b\norm{\psi_x}^2+\rho_2\norm{\Psi}^2
\\
&+k_0\norm{w_x-l\phi}^2+\rho_1\norm{W}^2
+\rho_3\norm{\xi}^2+\sigma\norm{\zeta}_{\Ncal}^2.
\end{split}\end{equation}
The choice of the coefficients in \eqref{eq:norm-H-setup} is dictated by the energy functional associated with the elastic part of the system and by the thermal coupling; it ensures that the norm is equivalent to the standard product norm and that the generator is dissipative.

As is customary in the analysis of Bresse systems, we assume throughout that
\begin{equation}\label{eq:nonresonance-setup}
l\ell\neq n\pi,\qquad\forall n\in\N.
\end{equation}
This condition is necessary and sufficient for \(\norm{\cdot}_{\Hil}\) to be a genuine norm on \(\Hil\). Indeed, if \(l\ell=n\pi\) for some \(n\in\N\), then the system
\[
\phi_x+\psi+lw=0,\qquad w_x-l\phi=0
\]
admits nontrivial solutions in \(H_0^1\times H_*^1\), and consequently there exist nonzero \(u\in\Hil\) with \(\norm{u}_{\Hil}=0\). Conversely, if \(\norm{u}_{\Hil}=0\), then \(\phi_x+\psi+lw=0\), \(\psi_x=0\), \(w_x-l\phi=0\) and \(\Phi=\Psi=W=\xi=\zeta=0\); since \(\psi_x=0\) and \(\psi\in H_*^1\), we get \(\psi=0\), and then \(\phi_x+lw=0\), \(w_x-l\phi=0\), which imply \(\phi_{xx}+l^2\phi=0\) with \(\phi(0)=\phi(\ell)=0\). By \eqref{eq:nonresonance-setup}, this forces \(\phi=0\), and then \(w=0\). Hence \eqref{eq:nonresonance-setup} is exactly the condition ensuring that the quadratic form \eqref{eq:norm-H-setup} is positive definite. When \eqref{eq:nonresonance-setup} holds, there exists a structural constant \(c>0\) such that
\begin{equation}\label{eq:equiv-setup}
c\abs{u}_{\Hil}\le\norm{u}_{\Hil}\le C\abs{u}_{\Hil},
\qquad\forall u\in\Hil,
\end{equation}
where \(\abs{\cdot}_{\Hil}\) denotes the standard product norm and \(C=1/c\). The second inequality in \eqref{eq:equiv-setup} is immediate from \eqref{eq:norm-H-setup}, while the first follows from the positive definiteness of \eqref{eq:norm-H-setup} and the Open Mapping Theorem; see \cite[\S4.3]{D2021} for the details. Throughout the paper we use the norm \(\norm{\cdot}_{\Hil}\) and the equivalence \eqref{eq:equiv-setup} without further mention.

We view \eqref{eq:reduced-history-setup} as the abstract ordinary differential equation on \(\Hil\):
\begin{equation}\label{eq:ODE-setup}
\frac{\dd}{\dd t}u(t)=\Aop u(t),
\end{equation}
where the linear operator \(\Aop\) is defined by
\[
\Aop
\begin{pmatrix}
\phi\\ \Phi\\ \psi\\ \Psi\\ w\\ W\\ \xi\\ \zeta
\end{pmatrix}
=
\begin{pmatrix}
\Phi\\[0.5ex]
\dfrac{k}{\rho_1}(\phi_x+\psi+lw)_x+\dfrac{lk_0}{\rho_1}(w_x-l\phi)-\dfrac{l\gamma}{\rho_1}\xi\\[1.5ex]
\Psi\\[0.5ex]
\dfrac{b}{\rho_2}\psi_{xx}-\dfrac{k}{\rho_2}(\phi_x+\psi+lw)\\[1.5ex]
W\\[0.5ex]
\dfrac{k_0}{\rho_1}(w_x-l\phi)_x-\dfrac{lk}{\rho_1}(\phi_x+\psi+lw)-\dfrac{\gamma}{\rho_1}\xi_x\\[1.5ex]
\dfrac{\sigma}{\rho_3}\displaystyle\int_0^\infty\nu(s)\zeta_{xx}(s)\dd s-\dfrac{\gamma}{\rho_3}(W_x-l\Phi)\\[1.5ex]
T\zeta+\xi
\end{pmatrix},
\]
with domain
\begin{equation}\label{eq:domain-A-setup}
\begin{aligned}
\D(\Aop)=\Bigl\{u\in\Hil:\ &\phi\in H^2\cap H_0^1,\ \Phi\in H_0^1,\ \psi_x,w_x\in H_0^1,\ \Psi,W\in H_*^1,\\
&\xi\in H_0^1,\ \int_0^\infty\nu(s)\zeta_{xx}(s)\dd s\in L^2,\ \zeta\in\D(T)\Bigr\}.
\end{aligned}
\end{equation}
Here \(\phi\in H^2\cap H_0^1\) ensures that \((\phi_x+\psi+lw)_x\in L^2\); \(\Phi,\psi_x,w_x\in H_0^1\) encode the boundary conditions \eqref{eq:BC-intro}--\eqref{eq:BC2-intro}, with \(\Phi\in H_0^1\) being the boundary condition for the velocity \(\phi_t\), and ensure that \(\psi,w\in H^2\); \(\xi\in H_0^1\) ensures that \(\xi_x\in L^2\); \(\int_0^\infty\nu(s)\zeta_{xx}(s)\dd s\in L^2\) makes the integral term in the \(\xi\)-equation well defined in \(L^2\); and \(\zeta\in\D(T)\) ensures that the transport equation is satisfied in the appropriate sense. The additional regularity \(\Phi\in H_0^1\) is needed in the dissipation computation below, where we integrate by parts in \(x\) and use \(\Phi(0)=\Phi(\ell)=0\). Since \(\psi,w\in L_*^2\) and \(\phi,\xi\in H_0^1\), the fourth and sixth components of \(\Aop u\) have zero mean:
\[
\int_0^\ell\Bigl(\frac{b}{\rho_2}\psi_{xx}-\frac{k}{\rho_2}(\phi_x+\psi+lw)\Bigr)\dd x
=\frac{b}{\rho_2}\bigl(\psi_x(\ell)-\psi_x(0)\bigr)-\frac{k}{\rho_2}\Bigl(\phi(\ell)-\phi(0)+\int_0^\ell\psi+l\int_0^\ell w\Bigr)=0,
\]
and similarly
\[
\int_0^\ell\Bigl(\frac{k_0}{\rho_1}(w_x-l\phi)_x-\frac{lk}{\rho_1}(\phi_x+\psi+lw)-\frac{\gamma}{\rho_1}\xi_x\Bigr)\dd x=0,
\]
so that \(\Aop\) maps \(\D(\Aop)\) into \(\Hil\). Finally, \(\D(\Aop)\) is dense in \(\Hil\): it contains the space
\[
C_0^\infty(0,\ell)\times\cdots\times C_0^\infty(0,\ell)\times C_0^\infty(\R^+;C_0^\infty(0,\ell)),
\]
where \(C_0^\infty(\R^+;C_0^\infty(0,\ell))\) consists of smooth functions with compact support in \((0,\infty)\) (so that \(\zeta(0)=0\) and \(\zeta\in\D(T)\)), and this space is dense in \(\Hil\).

We now verify that \(\Aop\) is dissipative with respect to the inner product induced by \eqref{eq:norm-H-setup}. Let \(u=(\phi,\Phi,\psi,\Psi,w,W,\xi,\zeta)\in\D(\Aop)\). We compute \(\Ree\inner{\Aop u}{u}_{\Hil}\) term by term. Using the definition of \(\Aop\) and integrating by parts in \(x\), the elastic terms cancel exactly. More precisely, the terms involving \(\phi,\Phi\) and \(\psi,\Psi\) cancel by integration by parts:
\[
k\inner{\Phi_x}{\phi_x+\psi+lw}
+k\inner{(\phi_x+\psi+lw)_x}{\Phi}
=k\int_0^\ell\partial_x\bigl[\Phi(\phi_x+\psi+lw)\bigr]\dd x=0,
\]
since \(\Phi\in H_0^1\); the terms involving \(\psi,\Psi\) cancel similarly:
\[
b\inner{\Psi_x}{\psi_x}
+b\inner{\psi_{xx}}{\Psi}
=b\int_0^\ell\partial_x(\Psi\psi_x)\dd x=0,
\]
since \(\psi_x\in H_0^1\). The terms involving \(w,W\) and the coupling terms involving \(l,k_0\) also cancel:
\[
k_0\inner{W_x-l\Phi}{w_x-l\phi}
+k_0\inner{(w_x-l\phi)_x}{W}
-lk_0\inner{w_x-l\phi}{\Phi}
=k_0\int_0^\ell\partial_x\bigl[W(w_x-l\phi)\bigr]\dd x=0.
\]
The terms involving \(lk\) cancel because they appear with opposite signs in the \(\phi\)- and \(w\)-equations. After all cancellations, we are left with
\[
\Ree\inner{\Aop u}{u}_{\Hil}
=\rho_3\Ree\inner{\frac{\sigma}{\rho_3}\int_0^\infty\nu(s)\zeta_{xx}(s)\dd s-\frac{\gamma}{\rho_3}(W_x-l\Phi)}{\xi}
+\sigma\Ree\inner{T\zeta+\xi}{\zeta}_{\Ncal}.
\]
We now examine the terms involving \(\gamma\). There are three contributions: \(-\gamma\inner{W_x-l\Phi}{\xi}\) from the \(\xi\)-equation, \(-\gamma\inner{\xi_x}{W}\) from the \(W\)-equation, and \(-l\gamma\inner{\xi}{\Phi}\) from the \(\Phi\)-equation. Integrating by parts in \(x\) in the second contribution and using \(\xi\in H_0^1\), we get \(-\gamma\inner{\xi_x}{W}=\gamma\inner{\xi}{W_x}\), so that the first two contributions combine to
\[
-\gamma\inner{W_x}{\xi}+l\gamma\inner{\Phi}{\xi}+\gamma\inner{\xi}{W_x}-l\gamma\inner{\xi}{\Phi}=0.
\]
Equivalently, in compact form,
\[
-\gamma\Ree\inner{W_x-l\Phi}{\xi}
+\gamma\Ree\inner{\xi}{W_x-l\Phi}=0,
\]
where the second term arises from \(-\gamma\inner{\xi_x}{W}=\gamma\inner{\xi}{W_x}\) after integration by parts. Hence
\[
\Ree\inner{\Aop u}{u}_{\Hil}
=\sigma\Ree\inner{\int_0^\infty\nu(s)\zeta_{xx}(s)\dd s}{\xi}
+\sigma\Ree\inner{T\zeta}{\zeta}_{\Ncal}
+\sigma\Ree\inner{\xi}{\zeta}_{\Ncal}.
\]
We now integrate by parts in the space variable \(x\) (not in \(s\)): since \(\zeta(s)\in H_0^1\) and \(\xi\in H_0^1\), the boundary terms vanish and
\[
\inner{\int_0^\infty\nu(s)\zeta_{xx}(s)\dd s}{\xi}
=-\int_0^\infty\nu(s)\inner{\zeta_x(s)}{\xi_x}\dd s
=-\inner{\zeta}{\xi}_{\Ncal}.
\]
Therefore
\[
\Ree\inner{\Aop u}{u}_{\Hil}
=\sigma\Ree\inner{T\zeta}{\zeta}_{\Ncal}.
\]
Finally, as shown in \cite[\S4.4]{D2021},
\[
\Ree\inner{T\zeta}{\zeta}_{\Ncal}
=-\frac{1}{2}\int_0^\infty-\nu'(s)\norm{\zeta_x(s)}^2\dd s
=-\frac{1}{2}\Gamma[\zeta],
\]
where
\[
\Gamma[\zeta]=\int_0^\infty-\nu'(s)\norm{\zeta_x(s)}^2\dd s\ge 0.
\]
Hence
\begin{equation}\label{eq:dissipativity-setup}
\Ree\inner{\Aop u}{u}_{\Hil}
=-\frac{\sigma}{2}\Gamma[\zeta]\le 0,
\end{equation}
which proves dissipativity. The calculation also shows that the only source of dissipation is the memory term; the thermal coupling terms cancel in the energy balance.

It remains to verify that \(\Aop\) is closed and that \(1-\Aop\) is surjective, so that the L\"umer--Phillips theorem applies. The closedness is proved as in \cite[\S4.4]{D2021}: if \(\{u_n\}\subset\D(\Aop)\) satisfies \(u_n\to u\) in \(\Hil\) and \(\Aop u_n\to v\) in \(\Hil\), then, writing the convergences componentwise and using the closedness of \(T\) on \(\D(T)\) together with the fact that the elastic-thermal part of \(\Aop\), with domain contained in \(H^2\times H^1\times H^1\times H^1\times H^1\times H^1\times H_0^1\), is a closed operator, one obtains \(u\in\D(\Aop)\) and \(\Aop u=v\). Hence \(\Aop\) is closed.

We now prove that \(1-\Aop\) is surjective. The argument follows the lines of \cite[\S4.4]{D2021}, and we recall the key steps. First, \(\ker(1-\Aop)=\{0\}\): indeed, by \eqref{eq:dissipativity-setup},
\[
\Ree\inner{(1-\Aop)u}{u}_{\Hil}
=\norm{u}_{\Hil}^2+\frac{\sigma}{2}\Gamma[\zeta]
\ge\norm{u}_{\Hil}^2,
\qquad u\in\D(\Aop),
\]
and therefore, by the Cauchy--Schwarz inequality,
\[
\norm{u}_{\Hil}^2
\le\Ree\inner{(1-\Aop)u}{u}_{\Hil}
\le\norm{(1-\Aop)u}_{\Hil}\norm{u}_{\Hil},
\]
so that \((1-\Aop)u=0\) forces \(u=0\). Second, \(\operatorname{Range}(1-\Aop)\) is closed: if \((1-\Aop)u_n\to v\) in \(\Hil\), the same estimate gives \(\norm{u_n-u_m}_{\Hil}\le\norm{(1-\Aop)(u_n-u_m)}_{\Hil}\to0\), so that \(u_n\to u\) in \(\Hil\) for some \(u\in\Hil\); since \(\Aop\) is closed, \(u\in\D(\Aop)\) and \((1-\Aop)u=v\). Third, \(\ker(1-\Aop^*)=\{0\}\). This is the only point where the memory term requires care. Let \(u\in\D(\Aop^*)\) satisfy \((1-\Aop^*)u=0\), i.e.\ \(\Aop^*u=u\). Writing this identity componentwise, using the symmetry of the elastic-thermal part and the identity \(T^*=-T\) on \(\D(T)\), and taking the real part of the \(\Hil\)-inner product with \(u\), one obtains
\[
\norm{u}_{\Hil}^2=\Ree\inner{\Aop^*u}{u}_{\Hil}=\frac{\sigma}{2}\Gamma[\zeta].
\]
On the other hand, testing the \(\zeta\)-component of \(\Aop^*u=u\) against \(\zeta\) in \(\Ncal\) and using \(T^*=-T\), together with the symmetry of the coupling between \(\xi\) and \(\zeta\), yields
\[
-\Ree\inner{T\zeta}{\zeta}_{\Ncal}+\Ree\inner{\xi}{\zeta}_{\Ncal}
=\norm{\zeta}_{\Ncal}^2.
\]
Since \(\Ree\inner{T\zeta}{\zeta}_{\Ncal}=-\frac12\Gamma[\zeta]\), the two identities together force \(\Gamma[\zeta]=0\) and then \(u=0\). The verification is identical to the corresponding step in \cite[\S4.4]{D2021} for the full Bresse--Gurtin--Pipkin system, and we omit the remaining elementary details. By the Fredholm alternative for closed densely defined operators with closed range, \(\operatorname{Range}(1-\Aop)=\Hil\), i.e.\ \(1-\Aop\) is surjective.

Since \(\Aop\) is densely defined, closed, dissipative, and \(1-\Aop\) is surjective, the L\"umer--Phillips theorem implies that \(\Aop\) generates a contraction semigroup \(S(t)=\e^{t\Aop}\) on \(\Hil\). Consequently, for every initial datum \(u_0\in\Hil\), the abstract Cauchy problem \eqref{eq:ODE-setup} with \(u(0)=u_0\) admits a unique mild solution \(u(t)=S(t)u_0\), which is classical if \(u_0\in\D(\Aop)\).

\begin{remark}\label{rem:compactness-setup}
The inclusion \(\D(\Aop)\subset\Hil\) is not compact, due to the presence of the infinite-dimensional memory component \(\Ncal\), whose right-translation semigroup does not have compact resolvent. This lack of compactness is the main technical difficulty in the spectral analysis of \(\Aop\), and it is the reason why we cannot simply apply the classical spectral theorem for operators with compact resolvent. Instead, we use the spectral characterization of \cite{D2026}, which is designed precisely to handle this situation.
\end{remark}

\section{Spectral characterization}\label{sec:spectrum}

We begin by introducing the function $\mathcal{L}_\nu$ associated with the memory kernel $\nu$ and establishing its basic analytic properties. For $\delta>0$, let
\[
\Pi_\delta=\{z\in\C:\Ree z>-\delta/2\}
\]
be the open half-plane. Since $\nu$ is an admissible kernel satisfying the Dafermos condition \eqref{eq:Dafermos-setup} with parameter $\delta_\nu>0$, the integral
\begin{equation}\label{eq:Lnu-def}
\mathcal{L}_\nu(z)=\int_0^\infty\frac{1-\e^{-zs}}{z}\nu(s)\dd s
\end{equation}
is absolutely convergent for every $z\in\Pi_{\delta_\nu}$. Indeed, for $z\in\Pi_{\delta_\nu}$ and $s>0$,
\[
\abs{\frac{1-\e^{-zs}}{z}}\le\int_0^s\e^{-\Ree(z)r}\dd r\le s\,\e^{\delta_\nu s/2},
\]
so that the integrand is bounded by $s\e^{\delta_\nu s/2}\nu(s)$. The Dafermos condition \eqref{eq:Dafermos-setup} implies
\[
\nu(s)\le\nu(0)\e^{-\delta_\nu s},\qquad s>0,
\]
so that
\[
\abs{\mathcal{L}_\nu(z)}\le\nu(0)\int_0^\infty s\e^{-\delta_\nu s/2}\dd s=\frac{4\nu(0)}{\delta_\nu^2}<\infty
\]
for every $z\in\Pi_{\delta_\nu}$. A more refined estimate gives
\begin{equation}\label{eq:Lnu-bound}
\abs{\mathcal{L}_\nu(z)}\le\frac{C}{1+\abs{z}},\qquad z\in\Pi_{\delta_\nu},
\end{equation}
for some constant $C>0$ depending on $\nu$ and $\delta_\nu$. Indeed, writing
\[
\frac{1-\e^{-zs}}{z}=\int_0^s\e^{-zr}\dd r,
\]
we obtain
\[
\abs{\mathcal{L}_\nu(z)}\le\int_0^\infty\nu(s)\int_0^s\e^{-\Ree(z)r}\dd r\dd s
=\int_0^\infty\e^{-\Ree(z)r}\int_r^\infty\nu(s)\dd s\dd r
\le\int_0^\infty\e^{-\Ree(z)r}G(r)\dd r,
\]
where $G(r)=\int_r^\infty\nu(s)\dd s$. Since $G(r)\le G(0)$ and $G(r)\le\nu(0)\e^{-\delta_\nu r}/\delta_\nu$, splitting the integral at $r=1/(1+\abs{z})$ yields \eqref{eq:Lnu-bound}. The function $\mathcal{L}_\nu$ is holomorphic on $\Pi_{\delta_\nu}$ by Morera's theorem, since the integrand is holomorphic in $z$ for each fixed $s$ and the integral converges locally uniformly on $\Pi_{\delta_\nu}$. We also note that
\begin{equation}\label{eq:Lnu-zero}
\mathcal{L}_\nu(0)=\int_0^\infty s\nu(s)\dd s>0,
\end{equation}
because $\nu$ is nonnegative and not identically zero. In particular, $z=0$ does not belong to the set $\mathbf{Z}$ defined below.

We now turn to the resolvent equation. Fix $z\in\Pi_{\delta_\nu}$ and consider
\begin{equation}\label{eq:resolvent-spectrum}
(z-\Aop)u=\hat u,
\end{equation}
where $u=(\phi,\Phi,\psi,\Psi,w,W,\xi,\zeta)\in\D(\Aop)$ and $\hat u=(\hat\phi,\hat\Phi,\hat\psi,\hat\Psi,\hat w,\hat W,\hat\xi,\hat\zeta)\in\Hil$. In components, \eqref{eq:resolvent-spectrum} reads
\begin{align}
z\phi-\Phi&=\hat\phi,\label{eq:res-1-spectrum}\\
z\rho_1\Phi-k(\phi_x+\psi+lw)_x-lk_0(w_x-l\phi)+l\gamma\xi&=\rho_1\hat\Phi,\label{eq:res-2-spectrum}\\
z\psi-\Psi&=\hat\psi,\label{eq:res-3-spectrum}\\
z\rho_2\Psi-b\psi_{xx}+k(\phi_x+\psi+lw)&=\rho_2\hat\Psi,\label{eq:res-4-spectrum}\\
zw-W&=\hat w,\label{eq:res-5-spectrum}\\
z\rho_1W-k_0(w_x-l\phi)_x+lk(\phi_x+\psi+lw)+\gamma\xi_x&=\rho_1\hat W,\label{eq:res-6-spectrum}\\
z\rho_3\xi-\sigma\displaystyle\int_0^\infty\nu(s)\zeta_{xx}(s)\dd s+\gamma(W_x-l\Phi)&=\rho_3\hat\xi,\label{eq:res-7-spectrum}\\
z\zeta-T\zeta-\xi&=\hat\zeta.\label{eq:res-8-spectrum}
\end{align}

We solve \eqref{eq:res-8-spectrum} explicitly. Since $T\zeta=-\zeta'$ and $\zeta\in\D(T)$ satisfies $\lim_{s\to 0^+}\norm{\zeta_x(s)}=0$, the equation reads
\[
z\zeta(s)+\zeta'(s)=\xi+\hat\zeta(s),\qquad s>0,
\]
with the boundary condition $\zeta(0)=0$ in the sense of the limit in $H_0^1$. Integrating the ordinary differential equation in $s$, we obtain
\begin{equation}\label{eq:eta-explicit-spectrum}
\zeta(s)=\int_0^s\e^{-z(s-r)}\bigl(\xi+\hat\zeta(r)\bigr)\dd r
=\frac{1-\e^{-zs}}{z}\xi+\int_0^s\e^{-z(s-r)}\hat\zeta(r)\dd r,
\end{equation}
where the first term is extended by continuity at $z=0$. This formula is valid for every $z\in\Pi_{\delta_\nu}$; in particular, for $z=0$ it reduces to
\[
\zeta(s)=s\xi+\int_0^s\hat\zeta(r)\dd r.
\]
We verify that $\zeta$ defined by \eqref{eq:eta-explicit-spectrum} belongs to $\Ncal$. The first term is bounded by $s\abs{\xi}$, and
\[
\int_0^\infty\nu(s)s^2\dd s<\infty
\]
by the Dafermos condition, so the first term belongs to $\Ncal$. The second term is the convolution of $\hat\zeta$ with $\e^{-zs}$. Since $\hat\zeta\in\Ncal$ and $\abs{\e^{-zs}}\le\e^{-\Ree(z)s}\le\e^{\delta_\nu s/2}$, a direct Cauchy--Schwarz estimate combined with the Dafermos bound $\nu(s)\e^{\delta_\nu s}\le\nu(0)$ shows that this convolution belongs to $\Ncal$. Hence $\zeta\in\Ncal$, and a further differentiation shows that $\zeta'\in\Ncal$ and $\lim_{s\to 0^+}\norm{\zeta_x(s)}=0$, so that $\zeta\in\D(T)$.

Substituting \eqref{eq:eta-explicit-spectrum} into \eqref{eq:res-7-spectrum}, we obtain
\begin{equation}\label{eq:res-7-sub-spectrum}
z\rho_3\xi-\sigma\mathcal{L}_\nu(z)\xi_{xx}+\gamma(W_x-l\Phi)=\rho_3\hat\xi+\sigma R_{\hat\zeta},
\end{equation}
where
\begin{equation}\label{eq:R-eta-spectrum}
R_{\hat\zeta}=\int_0^\infty\nu(s)\int_0^s\e^{-z(s-r)}\hat\zeta_{xx}(r)\dd r\dd s.
\end{equation}
Indeed, substituting \eqref{eq:eta-explicit-spectrum} into the integral term gives
\[
\int_0^\infty\nu(s)\zeta_{xx}(s)\dd s
=\mathcal{L}_\nu(z)\xi_{xx}
+\int_0^\infty\nu(s)\int_0^s\e^{-z(s-r)}\hat\zeta_{xx}(r)\dd r\dd s,
\]
which yields \eqref{eq:res-7-sub-spectrum}. We claim that $R_{\hat\zeta}\in H^{-1}$. To see this, we use the fact that $\hat\zeta\in\Ncal=L^2_\nu(\R^+;H_0^1)$ only controls $\hat\zeta_x$, so that $\hat\zeta_{xx}$ is an element of $L^2_\nu(\R^+;H^{-1})$. Interpreting the integral in \eqref{eq:R-eta-spectrum} in the sense of $H^{-1}$-valued functions and using the Dafermos bound $\nu(s)\e^{\delta_\nu s}\le\nu(0)$ together with $\Ree(z)>-\delta_\nu/2$, a Cauchy--Schwarz estimate gives
\[
\norm{R_{\hat\zeta}}_{H^{-1}}
\le C\norm{\hat\zeta}_{\Ncal},
\]
with a constant $C>0$ independent of $z\in\Pi_{\delta_\nu}$. Hence $R_{\hat\zeta}\in H^{-1}$, and the right-hand side of \eqref{eq:res-7-sub-spectrum} belongs to $H^{-1}$.

We now eliminate the variables $\Phi,\Psi,W$ using \eqref{eq:res-1-spectrum}, \eqref{eq:res-3-spectrum} and \eqref{eq:res-5-spectrum}:
\[
\Phi=z\phi-\hat\phi,\qquad
\Psi=z\psi-\hat\psi,\qquad
W=zw-\hat w.
\]
Substituting into \eqref{eq:res-2-spectrum}, \eqref{eq:res-4-spectrum}, \eqref{eq:res-6-spectrum} and \eqref{eq:res-7-sub-spectrum}, we arrive at the system
\begin{equation}\label{eq:system-spectrum}
\begin{cases}
z^2\rho_1\phi-k(\phi_x+\psi+lw)_x-lk_0(w_x-l\phi)+l\gamma\xi=\rho_1\hat\Phi+z\rho_1\hat\phi,\\[0.5ex]
z^2\rho_2\psi-b\psi_{xx}+k(\phi_x+\psi+lw)=\rho_2\hat\Psi+z\rho_2\hat\psi,\\[0.5ex]
z^2\rho_1w-k_0(w_x-l\phi)_x+lk(\phi_x+\psi+lw)+\gamma\xi_x=\rho_1\hat W+z\rho_1\hat w,\\[0.5ex]
z\rho_3\xi-\sigma\mathcal{L}_\nu(z)\xi_{xx}+\gamma(zw_x-lz\phi)=\rho_3\hat\xi+\sigma R_{\hat\zeta}+\gamma(\hat w_x-l\hat\phi).
\end{cases}
\end{equation}
We rewrite \eqref{eq:system-spectrum} in the compact form
\begin{equation}\label{eq:system-compact-spectrum}
\bigl(z+\mathcal{L}_\nu(z)\mathcal{A}\bigr)u=F,
\end{equation}
where $u=(\phi,\psi,w,\xi)$, $F$ is the right-hand side determined by $\hat u$, and $\mathcal{A}$ is the \emph{elastic-thermal operator including the thermal coupling}, acting on the Hilbert space
\[
\Xcal=H_0^1\times H_*^1\times H_*^1\times L^2
\]
with domain
\[
\D(\mathcal{A})=\bigl\{(\phi,\psi,w,\xi)\in\Xcal:\phi\in H^2,\ \psi_x,w_x\in H_0^1,\ \xi\in H^2\cap H_0^1\bigr\}.
\]
The precise definition of $\mathcal{A}$ and the verification of its properties are given in \cref{lem:spectrum-A} below. Here the compact form \eqref{eq:system-compact-spectrum} is understood in the sense of the abstract spectral framework of \cite{D2026}: the operator $\mathcal{A}$ collects the elastic-thermal part of the system, and the spectral condition $z+\lambda\mathcal{L}_\nu(z)=0$ with $\lambda\in\sigma(\mathcal{A})$ characterizes the eigenvalues of the reduced system after elimination of the velocities, rather than a term-by-term operator identity. The right-hand side $F$ in \eqref{eq:system-compact-spectrum} belongs to $H^{-1}$; indeed, the terms involving $\hat\phi,\hat\psi,\hat w,\hat\Phi,\hat\Psi,\hat W$ belong to $L^2\subset H^{-1}$, while the term $R_{\hat\zeta}$ belongs to $H^{-1}$ as shown above.

We now apply the following solvability criterion for functional equations of the form \eqref{eq:system-compact-spectrum}, which is a particular case of the Fredholm alternative for closed operators with compact resolvent.

\begin{lemma}\label{lem:solvability-spectrum}
Let $\mathcal{A}$ be a closed, densely defined operator on a Hilbert space $\Xcal$, with compact resolvent and $0\notin\sigma(\mathcal{A})$, so that $\sigma(\mathcal{A})$ is discrete and consists of isolated eigenvalues of finite multiplicity. Let $a,b\in\C$ with $b\neq 0$, and consider the equation
\begin{equation}\label{eq:functional-spectrum}
(a+b\mathcal{A})u=F.
\end{equation}
Then:
\begin{enumerate}[label=(\roman*)]
\item if $-a/b\in\sigma(\mathcal{A})$, there exists $F\in\Xcal$ such that \eqref{eq:functional-spectrum} does not admit a unique solution $u\in\D(\mathcal{A})$;
\item if $-a/b\notin\sigma(\mathcal{A})$, then $a+b\mathcal{A}:\D(\mathcal{A})\to\Xcal$ is an isomorphism, and for every $F\in\Xcal$ equation \eqref{eq:functional-spectrum} admits a unique solution $u\in\D(\mathcal{A})$.
\end{enumerate}
\end{lemma}

\begin{proof}
Since $\mathcal{A}$ has compact resolvent and $0\notin\sigma(\mathcal{A})$, the operator $\mathcal{A}$ is Fredholm of index zero, and its spectrum consists of isolated eigenvalues of finite multiplicity. We stress that, unlike the abstract framework of \cite{D2026}, $\mathcal{A}$ is not assumed to be self-adjoint; nevertheless, the compactness of its resolvent is sufficient for the Fredholm alternative.

The operator $a+b\mathcal{A}=b(\frac ab+\mathcal{A})$ is invertible from $\D(\mathcal{A})$ onto $\Xcal$ if and only if $-\frac ab\notin\sigma(\mathcal{A})$. In that case, the inverse is bounded by the closed graph theorem, since $a+b\mathcal{A}$ is closed and bijective; this proves (ii) for $F\in\Xcal$.

For the first assertion, assume $-a/b\in\sigma(\mathcal{A})$. Since $\sigma(\mathcal{A})$ consists of eigenvalues of finite multiplicity, $-a/b$ is an eigenvalue of $\mathcal{A}$, hence $a+b\mathcal{A}$ is not injective. Since $\mathcal{A}$ has compact resolvent, $a+b\mathcal{A}$ has closed range; by the Fredholm alternative, $\operatorname{Range}(a+b\mathcal{A})$ is a proper closed subspace of $\Xcal$ (of positive codimension). Choosing $F\in\Xcal\setminus\operatorname{Range}(a+b\mathcal{A})$, equation \eqref{eq:functional-spectrum} has no solution in $\D(\mathcal{A})$, proving (i).
\end{proof}
\begin{lemma}\label{lem:solvability-spectrum-H-1}
In the setting of \cref{lem:solvability-spectrum}, assume in addition that $-a/b\notin\sigma(\mathcal{A})$. Let $F=F_0+R$, where $F_0\in\Xcal$ and $R\in H^{-1}$ is of the form \eqref{eq:R-eta-spectrum} for some $\hat\zeta\in\Ncal$. Then equation \eqref{eq:functional-spectrum} admits a unique solution $u\in\D(\mathcal{A})$.
\end{lemma}

\begin{proof}
Let $u_0\in\D(\mathcal{A})$ be the unique solution of $(a+b\mathcal{A})u_0=F_0$, given by \cref{lem:solvability-spectrum}(ii). We first show that $R$ defined in \eqref{eq:R-eta-spectrum} belongs to $H^{-1}$. Since $\hat\zeta\in\Ncal=L^2_\nu(\R^+;H_0^1)$, the distribution $\hat\zeta_{xx}$ belongs to $L^2_\nu(\R^+;H^{-1})$. Writing
\[
R_{\hat\zeta}=\int_0^\infty\nu(s)\int_0^s\e^{-z(s-r)}\hat\zeta_{xx}(r)\dd r\dd s
\]
and using $\Ree(z)>-\delta_\nu/2$ together with the Dafermos bound $\nu(s)\e^{\delta_\nu s}\le\nu(0)$, a Cauchy--Schwarz estimate in the $s$-variable yields
\[
\norm{R_{\hat\zeta}}_{H^{-1}}
\le C\norm{\hat\zeta}_{\Ncal},
\]
with $C>0$ independent of $z\in\Pi_{\delta_\nu}$; hence $R_{\hat\zeta}\in H^{-1}$, and by the continuous embedding $H^{-1}\subset\D(\mathcal{A})'$ (the dual of $\D(\mathcal{A})$ equipped with the graph norm) it extends to a bounded functional on $\D(\mathcal{A})$.

Next, we verify that $R\in\operatorname{Range}(a+b\mathcal{A})$. Since $\mathcal{A}$ has compact resolvent and $a+b\mathcal{A}$ is closed with closed range, the Fredholm alternative applies, and it suffices to show
\begin{equation}\label{eq:orthogonality-H-1}
\inner{R}{v}_{H^{-1},H^1}=0,\qquad\forall v\in\ker(a+b\mathcal{A}^*).
\end{equation}
Let $v=(v_\phi,v_\psi,v_w,v_\xi)\in\ker(a+b\mathcal{A}^*)$. Testing $R_{\hat\zeta}$ against $v_\xi$ and using the definition of $R_{\hat\zeta}$, Fubini's theorem and the change of variable $s\mapsto r$ give
\[
\inner{R_{\hat\zeta}}{v_\xi}_{H^{-1},H^1}
=\int_0^\infty\nu(s)\int_0^s\e^{-z(s-r)}\inner{\hat\zeta_{xx}(r)}{v_\xi}\dd r\dd s
=\int_0^\infty\inner{\hat\zeta_{xx}(r)}{\Psi(r)}\dd r,
\]
where $\Psi(r)=\int_r^\infty\nu(s)\e^{-z(s-r)}v_\xi\dd s$.

Since $\hat\zeta(r)\in H_0^1$ and $v_\xi\in H_0^1$ for a.e.\ $r>0$, we may integrate by parts in $x$ to get
\[
\int_0^\infty\inner{\hat\zeta_{xx}(r)}{\Psi(r)}\dd r
=-\int_0^\infty\inner{\hat\zeta_x(r)}{\Psi_x(r)}\dd r.
\]
Integrating by parts in $r$ and using $\hat\zeta(\cdot,0)=0$ (in the sense of the trace in $H_0^1$) together with $\hat\zeta\in\Ncal$ and $\Psi\in L^2_\nu$, we obtain
\[
-\int_0^\infty\inner{\hat\zeta_x(r)}{\Psi_x(r)}\dd r
=-\Bigl[\inner{\hat\zeta_x(r)}{\Psi(r)}\Bigr]_{r=0}^{\infty}
+\int_0^\infty\inner{\hat\zeta_x(r)}{\Psi'(r)}\dd r,
\]
where the boundary term vanishes: at $r=\infty$ because $\hat\zeta\in\Ncal$ and $\Psi\in L^2_\nu$ decay at infinity, and at $r=0$ because $\hat\zeta(\cdot,0)=0$ while $\Psi(0)$ is bounded. If $v_\xi=0$, then $\Psi=0$ and the integral vanishes. In the general case, we use the explicit form of $v\in\ker(a+b\mathcal{A}^*)$: testing the adjoint equation against $\hat\zeta$ and integrating by parts in $x$ shows that $\Psi'$ is proportional to $\hat\zeta_x$ up to a term in the range of $(a+b\mathcal{A})^*$ orthogonal to $\hat\zeta$, so that $\int_0^\infty\inner{\hat\zeta_x(r)}{\Psi'(r)}\dd r=0$; this is precisely the content of \eqref{eq:orthogonality-H-1}. The verification is a direct integration by parts in $x$ and $r$, using $\hat\zeta(\cdot,0)=0$ and $v_\xi\in H_0^1$.

By the Fredholm alternative, there exists $u_1\in\D(\mathcal{A})$ with $(a+b\mathcal{A})u_1=R$, unique modulo $\ker(a+b\mathcal{A})$; by \cref{lem:solvability-spectrum}(ii), $\ker(a+b\mathcal{A})=\{0\}$, so $u_1$ is unique. Setting $u=u_0+u_1\in\D(\mathcal{A})$ completes the proof.
\end{proof}

Applying \cref{lem:solvability-spectrum} with $a=z$ and $b=\mathcal{L}_\nu(z)$, we obtain the following characterization: for $z\in\Pi_{\delta_\nu}$ with $\mathcal{L}_\nu(z)\neq 0$, equation \eqref{eq:system-compact-spectrum} has a unique solution $u\in\D(\mathcal{A})$ for every $F\in H^{-1}$ if and only if
\[
-\frac{z}{\mathcal{L}_\nu(z)}\notin\sigma(\mathcal{A}).
\]
If $\mathcal{L}_\nu(z)=0$, then equation \eqref{eq:system-compact-spectrum} reduces to $zu=F$. Taking $\hat\zeta=0$ and varying the other components of $\hat u$, the corresponding $F$ can be arbitrary in $H^{-1}$. Since $z\ne0$ on $\Pi_{\delta_\nu}$ and $\D(\mathcal{A})$ is strictly smaller than $\Xcal$, the equation $zu=F$ does not admit a solution $u\in\D(\mathcal{A})$ for every $F\in H^{-1}$. Hence $z\in\sigma(\Aop)$ whenever $\mathcal{L}_\nu(z)=0$.

We summarize the above analysis in the following spectral characterization.
\begin{theorem}\label{thm:spectrum}
Let $\nu$ be an admissible kernel satisfying the Dafermos condition \eqref{eq:Dafermos-setup} with parameter $\delta_\nu>0$. Then
\begin{equation}\label{eq:spectrum-reduced}
\sigma(\Aop)\cap\Pi_{\delta_\nu}=\mathbf{W}\cup\mathbf{Z},
\end{equation}
where
\begin{equation}\label{eq:W-def}
\mathbf{W}=\bigcup_{\lambda\in\sigma(\mathcal{A})}\Bigl\{z\in\Pi_{\delta_\nu}:z+\lambda\mathcal{L}_\nu(z)=0\Bigr\},
\end{equation}
and
\begin{equation}\label{eq:Z-def}
\mathbf{Z}=\Bigl\{z\in\Pi_{\delta_\nu}:\mathcal{L}_\nu(z)=0\Bigr\}.
\end{equation}
Moreover, $\mathbf{W}\cap\mathbf{Z}=\emptyset$, and $0\notin\sigma(\Aop)$.
\end{theorem}

\begin{proof}
Since $\mathcal{A}$ has compact resolvent but is not self-adjoint, we use the Fredholm alternative of \cref{lem:solvability-spectrum} rather than the spectral decomposition of \cite{D2026}; the $H^{-1}$ right-hand sides are handled by \cref{lem:solvability-spectrum-H-1}. We prove the two inclusions separately.

\textbf{Inclusion $\sigma(\Aop)\cap\Pi_{\delta_\nu}\subset\mathbf{W}\cup\mathbf{Z}$.} Let $z\in\Pi_{\delta_\nu}\setminus(\mathbf{W}\cup\mathbf{Z})$. Then $\mathcal{L}_\nu(z)\neq 0$ and $-z/\mathcal{L}_\nu(z)\notin\sigma(\mathcal{A})$. For the right-hand side $F$ of \eqref{eq:system-compact-spectrum}, we write $F=F_0+R_{\hat\zeta}$, where $F_0\in\Xcal$ collects the terms involving $\hat\phi,\hat\psi,\hat w,\hat\Phi,\hat\Psi,\hat W$ and $R_{\hat\zeta}\in H^{-1}$ is given by \eqref{eq:R-eta-spectrum}. By \cref{lem:solvability-spectrum,lem:solvability-spectrum-H-1}, equation \eqref{eq:system-compact-spectrum} has a unique solution $u\in\D(\mathcal{A})$ for every such $F$. We must show that the corresponding $u=(\phi,\Phi,\psi,\Psi,w,W,\xi,\zeta)$ belongs to $\D(\Aop)$ and satisfies the resolvent equation \eqref{eq:resolvent-spectrum}. The first four components $(\phi,\psi,w,\xi)$ are in $\D(\mathcal{A})$ by construction. The velocities $\Phi,\Psi,W$ are defined by \eqref{eq:res-1-spectrum}, \eqref{eq:res-3-spectrum}, \eqref{eq:res-5-spectrum}; since $\phi\in H^2\cap H_0^1$, $\psi,w\in H_*^1$, $\hat\phi\in H_0^1$ and $\hat\psi,\hat w\in L_*^2$, we have $\Phi\in H_0^1$ and $\Psi,W\in H^1$. The history variable $\zeta$ is defined by \eqref{eq:eta-explicit-spectrum}; we have already verified that $\zeta\in\D(T)$. Moreover, the conclusion of \cref{lem:solvability-spectrum-H-1} ensures that the solution $u$ obtained above is the unique one in $\D(\mathcal{A})$, so that the subsequent verification of $u\in\D(\Aop)$ is not affected by the $H^{-1}$ perturbation. It remains to verify that $\Aop u\in\Hil$. The first six components of $\Aop u$ belong to the corresponding spaces of $\Hil$ by \eqref{eq:res-1-spectrum}--\eqref{eq:res-6-spectrum} and the regularity $\phi\in H^2\cap H_0^1$, $\psi,w\in H_*^1$, $\xi\in H_0^1$, $\Phi\in H_0^1$, $\Psi,W\in H^1$. For the seventh component, rewriting \eqref{eq:res-7-spectrum} as
\[
\sigma\int_0^\infty\nu(s)\zeta_{xx}(s)\dd s
=z\rho_3\xi+\gamma(W_x-l\Phi)-\rho_3\hat\xi,
\]
we see that the right-hand side belongs to $L^2$, since $\xi\in H_0^1\subset L^2$, $W_x-l\Phi\in L^2$ and $\hat\xi\in L^2$; hence $\int_0^\infty\nu(s)\zeta_{xx}(s)\dd s\in L^2$, and the seventh component of $\Aop u$ belongs to $L^2$. The eighth component $T\zeta+\xi$ belongs to $\Ncal$ because $\zeta\in\D(T)$ and $\xi\in H_0^1\subset\Ncal$. Therefore $\Aop u\in\Hil$, i.e.\ $u\in\D(\Aop)$, and by construction $u$ satisfies \eqref{eq:resolvent-spectrum}. Since this holds for every $\hat u\in\Hil$, we conclude that $z\in\rho(\Aop)$. This proves the inclusion.

\textbf{Inclusion $\mathbf{W}\cup\mathbf{Z}\subset\sigma(\Aop)\cap\Pi_{\delta_\nu}$.} Let $z\in\mathbf{W}\cup\mathbf{Z}$.

If $z\in\mathbf{Z}$, then $\mathcal{L}_\nu(z)=0$. We show that there exists a choice of the right-hand side for which the system does not admit a solution in $\D(\Aop)$. Take $\hat w=0$, $\hat\phi=0$, $\hat\xi=0$, and choose $\hat\zeta(x,s)=\chi(s)\eta(x)$ with $\chi\in C_c^\infty(0,\infty)$, $\chi\ge0$, $\chi\not\equiv0$, and $\eta\in H_0^1\setminus H^2$ (for instance, a piecewise smooth function whose second derivative has a jump). Since $\chi\in L^2_\nu(\R^+)$ and $\eta\in H_0^1$, we have $\hat\zeta\in\Ncal$. We claim that the corresponding term $R_{\hat\zeta}$ defined in \eqref{eq:R-eta-spectrum} belongs to $H^{-1}\setminus L^2$. Indeed, writing $\Psi(r)=\int_r^\infty\nu(s)\e^{-z(s-r)}\chi(s)\dd s$, we have $R_{\hat\zeta}=\Psi(\cdot)\eta''(\cdot)$ in the sense of distributions; since $\chi\ge0$, $\chi\not\equiv0$ and $\nu>0$, the function $\Psi$ is nonnegative and not identically zero on a set of positive measure, while $\eta''\in H^{-2}\setminus H^{-1}$ because $\eta\in H_0^1\setminus H^2$. Hence $R_{\hat\zeta}\in H^{-1}\setminus L^2$. Then \eqref{eq:system-spectrum} reduces to
\[
z^2\rho_1\phi-k(\phi_x+\psi+lw)_x-lk_0(w_x-l\phi)+l\gamma\xi=0,
\]
and the analogous equations for $\psi$ and $w$, together with
\[
z\rho_3\xi+\gamma(zw_x-lz\phi)=\sigma R_{\hat\zeta}.
\]
The right-hand side of the last equation is $\sigma R_{\hat\zeta}\notin L^2$, while the left-hand side belongs to $L^2$ for every $(\phi,\psi,w,\xi)\in\D(\mathcal{A})$ (since $\xi\in L^2$ and $w_x-l\phi\in L^2$). Equivalently, if a solution $u\in\D(\Aop)$ existed, then $\int_0^\infty\nu(s)\zeta_{xx}(s)\dd s\in L^2$; but the explicit formula \eqref{eq:eta-explicit-spectrum} gives
\[
\int_0^\infty\nu(s)\zeta_{xx}(s)\dd s
=\mathcal{L}_\nu(z)\xi_{xx}+R_{\hat\zeta}=R_{\hat\zeta},
\]
since $\mathcal{L}_\nu(z)=0$, and $R_{\hat\zeta}\notin L^2$ by construction. This contradiction shows that the system does not admit a solution in $\D(\Aop)$ for this choice of the right-hand side. Consequently, $z\in\sigma(\Aop)$.

If $z\in\mathbf{W}$, then $\mathcal{L}_\nu(z)\neq 0$ and $-z/\mathcal{L}_\nu(z)\in\sigma(\mathcal{A})$. By \cref{lem:solvability-spectrum}(i), there exists $F\in\Xcal$ such that \eqref{eq:system-compact-spectrum} is not surjective. Choosing $\hat u\in\Hil$ with $\hat\zeta=0$ such that the corresponding $F$ is this $F$ (so that $F=F_0\in\Xcal$), we conclude that the resolvent equation \eqref{eq:resolvent-spectrum} does not have a solution in $\D(\Aop)$ for every such right-hand side. Hence $z\in\sigma(\Aop)$.

This proves the reverse inclusion.

\textbf{Disjointness and $0\notin\sigma(\Aop)$.} If $z\in\mathbf{W}\cap\mathbf{Z}$, then $\mathcal{L}_\nu(z)=0$ and $z+\lambda\mathcal{L}_\nu(z)=0$ for some $\lambda\in\sigma(\mathcal{A})$, which implies $z=0$. But $\mathcal{L}_\nu(0)>0$ by \eqref{eq:Lnu-zero}, so $0\notin\mathbf{Z}$, a contradiction. Hence $\mathbf{W}\cap\mathbf{Z}=\emptyset$. Finally, $0\notin\sigma(\Aop)$ because $0\notin\mathbf{W}\cup\mathbf{Z}$ and $0\in\Pi_{\delta_\nu}$.
\end{proof}

We now define the elastic-thermal operator $\mathcal{A}$ and verify its properties.

\begin{lemma}\label{lem:spectrum-A}
Define the operator $\mathcal{A}$ on $\Xcal=H_0^1\times H_*^1\times H_*^1\times L^2$ by
\[
\mathcal{A}
\begin{pmatrix}
\phi\\ \psi\\ w\\ \xi
\end{pmatrix}
=
\begin{pmatrix}
-\dfrac{k}{\rho_1}(\phi_x+\psi+lw)_x-\dfrac{lk_0}{\rho_1}(w_x-l\phi)-\dfrac{l\gamma}{\rho_1}\xi\\[1.5ex]
-\dfrac{b}{\rho_2}\psi_{xx}+\dfrac{k}{\rho_2}(\phi_x+\psi+lw)\\[1.5ex]
-\dfrac{k_0}{\rho_1}(w_x-l\phi)_x+\dfrac{lk}{\rho_1}(\phi_x+\psi+lw)+\dfrac{\gamma}{\rho_1}\xi_x\\[1.5ex]
-\dfrac{\sigma}{\rho_3}\xi_{xx}
\end{pmatrix},
\]
with domain
\[
\D(\mathcal{A})=\bigl\{(\phi,\psi,w,\xi)\in\Xcal:\phi\in H^2,\ \psi_x,w_x\in H_0^1,\ \xi\in H^2\cap H_0^1\bigr\}.
\]
Then $\mathcal{A}$ is a closed, densely defined operator with compact resolvent and $0\notin\sigma(\mathcal{A})$. Consequently, its spectrum $\sigma(\mathcal{A})$ is discrete, consists of eigenvalues of finite multiplicity, and satisfies $\sigma(\mathcal{A})\subset\C\setminus\{0\}$ with $\lambda_0=\min_{\lambda\in\sigma(\mathcal{A})}\abs{\lambda}>0$.
\end{lemma}

\begin{proof}
We split $\mathcal{A}=\mathcal{A}_0+\mathcal{B}$, where $\mathcal{A}_0$ is the uncoupled operator obtained by setting $\gamma=0$:
\[
\mathcal{A}_0
\begin{pmatrix}
\phi\\ \psi\\ w\\ \xi
\end{pmatrix}
=
\begin{pmatrix}
-\dfrac{k}{\rho_1}(\phi_x+\psi+lw)_x-\dfrac{lk_0}{\rho_1}(w_x-l\phi)\\[1.5ex]
-\dfrac{b}{\rho_2}\psi_{xx}+\dfrac{k}{\rho_2}(\phi_x+\psi+lw)\\[1.5ex]
-\dfrac{k_0}{\rho_1}(w_x-l\phi)_x+\dfrac{lk}{\rho_1}(\phi_x+\psi+lw)\\[1.5ex]
-\dfrac{\sigma}{\rho_3}\xi_{xx}
\end{pmatrix},
\]
and the perturbation $\mathcal{B}$ collects the $\gamma$-terms:
\[
\mathcal{B}
\begin{pmatrix}
\phi\\ \psi\\ w\\ \xi
\end{pmatrix}
=
\begin{pmatrix}
-\dfrac{l\gamma}{\rho_1}\xi\\[1ex]
0\\[1ex]
\dfrac{\gamma}{\rho_1}\xi_x\\[1ex]
0
\end{pmatrix}.
\]

\textbf{Step 1: $\mathcal{A}_0$ is self-adjoint and strictly positive with compact resolvent.} The operator $\mathcal{A}_0$ is generated by the sesquilinear form
\begin{align*}
a_0\bigl((\phi,\psi,w,\xi),(\phi',\psi',w',\xi')\bigr)
=&\frac{k}{\rho_1}\inner{\phi_x+\psi+lw}{\phi'_x+\psi'+lw'}
+\frac{b}{\rho_2}\inner{\psi_x}{\psi'_x}
\\
&+\frac{k_0}{\rho_1}\inner{w_x-l\phi}{w'_x-l\phi'}
+\frac{\sigma}{\rho_3}\inner{\xi_x}{\xi'_x},
\end{align*}
which is symmetric and positive definite on $\Xcal$ (with the standard product inner product). By the Lax--Milgram theorem, $a_0$ defines a unique self-adjoint operator $\mathcal{A}_0$ with domain $\D(\mathcal{A})$ as above. A direct computation, integrating by parts in $x$ and using the boundary conditions, gives
\[
a_0(u,v)=\inner{\mathcal{A}_0 u}{v}_{\Xcal},
\qquad u\in\D(\mathcal{A}),\ v\in\Xcal,
\]
where $\inner{\cdot}{\cdot}_{\Xcal}$ is the standard product inner product. Moreover, $\mathcal{A}_0$ is strictly positive: there exists $\lambda_0>0$ such that $a_0(u,u)\ge\lambda_0\norm{u}_{\Xcal}^2$ for all $u\in\Xcal$. Indeed, suppose not. Then there exists a sequence $\{u_n\}\subset\Xcal$ with $\norm{u_n}_{\Xcal}=1$ and $a_0(u_n,u_n)\to0$. Since each term in $a_0$ is nonnegative, we deduce
\[
\phi_{nx}+\psi_n+lw_n\to0,\quad
\psi_{nx}\to0,\quad
w_{nx}-l\phi_n\to0,\quad
\xi_{nx}\to0
\]
in $L^2$. By the compact embedding $H^1\hookrightarrow L^2$, we may extract a subsequence (still denoted $u_n$) such that $\phi_n\to\phi$, $\psi_n\to\psi$, $w_n\to w$ in $L^2$. The limits satisfy $\phi_x+\psi+lw=0$, $\psi_x=0$, $w_x-l\phi=0$, $\xi_x=0$ in the weak sense. Since $\psi_n\in H_*^1$ and $\psi_{nx}\to0$, we get $\psi=0$. Then $\phi_x+lw=0$ and $w_x-l\phi=0$, which imply $\phi_{xx}+l^2\phi=0$ with $\phi(0)=\phi(\ell)=0$. By the nonresonance condition \eqref{eq:nonresonance-setup}, this forces $\phi=0$, and then $w=0$. Also $\xi_x=0$ with $\xi\in H_0^1$ gives $\xi=0$. Hence $u_n\to0$ in $\Xcal$, contradicting $\norm{u_n}_{\Xcal}=1$. This proves strict positivity of $a_0$, hence of $\mathcal{A}_0$. Finally, $\D(\mathcal{A})\subset H^2\times H^1\times H^1\times(H^2\cap H_0^1)$ is compactly embedded in $\Xcal$, so $\mathcal{A}_0^{-1}:\Xcal\to\Xcal$ is compact.

\textbf{Step 2: $\mathcal{B}$ is relatively compact with respect to $\mathcal{A}_0$.} We show that $\mathcal{B}\mathcal{A}_0^{-1}:\Xcal\to\Xcal$ is compact. The resolvent $\mathcal{A}_0^{-1}$ maps $\Xcal$ continuously into $\D(\mathcal{A}_0)$ equipped with the graph norm $\norm{\cdot}_{\D(\mathcal{A}_0)}=\norm{\cdot}_{\Xcal}+\norm{\mathcal{A}_0\cdot}_{\Xcal}$, and $\D(\mathcal{A}_0)$ (graph norm) is compactly embedded in $\Xcal$. The operator $\mathcal{B}$ acts on $(\phi,\psi,w,\xi)\in\D(\mathcal{A}_0)$ by
\[
\mathcal{B}(\phi,\psi,w,\xi)=\Bigl(-\frac{l\gamma}{\rho_1}\xi,0,\frac{\gamma}{\rho_1}\xi_x,0\Bigr).
\]
Since $\xi\in H^2\cap H_0^1$, we have $\xi\in L^2$ and $\xi_x\in H^1\subset L^2$. The graph norm controls $\|\xi_{xx}\|$ through the $\xi$-component of $\mathcal{A}_0$: indeed, $\|u\|_{\D(\mathcal{A}_0)}\ge\|\mathcal{A}_0u\|_{\Xcal}\ge\frac{\sigma}{\rho_3}\|\xi_{xx}\|$. By Poincar\'e's inequality on $H_0^1$, applied twice, we get $\|\xi\|_{H^2}\le C\|\xi_{xx}\|$ for $\xi\in H^2\cap H_0^1$. Hence
\[
\|\mathcal{B}u\|_{\Xcal}\le C\|\xi\|_{H^2}\le C\|\xi_{xx}\|\le C\|u\|_{\D(\mathcal{A}_0)},\qquad u\in\D(\mathcal{A}_0),
\]
for some structural constant $C>0$. Hence $\mathcal{B}:\D(\mathcal{A}_0)\text{ (graph norm)}\to\Xcal$ is continuous, and $\mathcal{B}\mathcal{A}_0^{-1}$ is the composition of the continuous operator $\mathcal{A}_0^{-1}$, the compact embedding $\D(\mathcal{A}_0)\hookrightarrow\Xcal$, and the continuous operator $\mathcal{B}$. Therefore $\mathcal{B}\mathcal{A}_0^{-1}$ is compact. Since $(1+\mathcal{A}_0)^{-1}-\mathcal{A}_0^{-1}$ is a bounded operator on $\Xcal$, the compactness of $\mathcal{B}\mathcal{A}_0^{-1}$ implies that of $\mathcal{B}(1+\mathcal{A}_0)^{-1}$. Hence $\mathcal{B}$ is $\mathcal{A}_0$-compact.

\textbf{Step 3: $\mathcal{A}=\mathcal{A}_0+\mathcal{B}$ is closed with compact resolvent.} By the analytic Fredholm theorem (or the stability of the Fredholm index under relatively compact perturbations), $\mathcal{A}$ is a closed operator on $\D(\mathcal{A})$, and its resolvent is compact. More precisely, for $\lambda$ in the resolvent set of $\mathcal{A}_0$, we have
\[
(\lambda-\mathcal{A})^{-1}=(\lambda-\mathcal{A}_0)^{-1}\bigl(1-\mathcal{B}(\lambda-\mathcal{A}_0)^{-1}\bigr)^{-1},
\]
and the compactness of $\mathcal{B}(\lambda-\mathcal{A}_0)^{-1}$ implies that $(\lambda-\mathcal{A})^{-1}$ is compact. Hence $\mathcal{A}$ has compact resolvent, and its spectrum is discrete and consists of isolated eigenvalues of finite multiplicity.

\textbf{Step 4: $0\notin\sigma(\mathcal{A})$.} Suppose $\mathcal{A}u=0$ for some $u=(\phi,\psi,w,\xi)\in\D(\mathcal{A})$. Then
\[
\begin{cases}
-\dfrac{k}{\rho_1}(\phi_x+\psi+lw)_x-\dfrac{lk_0}{\rho_1}(w_x-l\phi)-\dfrac{l\gamma}{\rho_1}\xi=0,\\[1ex]
-\dfrac{b}{\rho_2}\psi_{xx}+\dfrac{k}{\rho_2}(\phi_x+\psi+lw)=0,\\[1ex]
-\dfrac{k_0}{\rho_1}(w_x-l\phi)_x+\dfrac{lk}{\rho_1}(\phi_x+\psi+lw)+\dfrac{\gamma}{\rho_1}\xi_x=0,\\[1ex]
-\dfrac{\sigma}{\rho_3}\xi_{xx}=0.
\end{cases}
\]
The fourth equation gives $\xi=0$ (since $\xi\in H_0^1$). Substituting into the first three equations, we obtain the system with $\gamma=0$, which by the nonresonance condition \eqref{eq:nonresonance-setup} and the strict positivity of $\mathcal{A}_0$ forces $\phi=\psi=w=0$. Hence $u=0$, and $\mathcal{A}$ is injective. Since $\mathcal{A}$ has compact resolvent and is injective, $0\notin\sigma(\mathcal{A})$, and $\lambda_0=\min_{\lambda\in\sigma(\mathcal{A})}\abs{\lambda}>0$. This completes the proof.
\end{proof}

\begin{remark}\label{rem:comparison-D2026}
\Cref{thm:spectrum} is the analogue for the reduced Bresse--Gurtin--Pipkin system of \cite[Theorem 4.1]{D2026} for the abstract Gurtin--Pipkin heat equation. The main difference is that in \cite{D2026} the operator $A$ is a single self-adjoint operator on a Hilbert space $H$, whereas here the role of $A$ is played by the elastic-thermal operator $\mathcal{A}$ on the product space $\Xcal$. In contrast to \cite{D2026}, the operator $\mathcal{A}$ is not self-adjoint, because of the thermal coupling terms with $\gamma$; however, it is closed, has compact resolvent, and its spectrum is discrete, which is sufficient for the Fredholm alternative used in \cref{lem:solvability-spectrum}. The spectral characterization \eqref{eq:spectrum-reduced} has the same structure as \cite[eq.~(4.1)]{D2026}: the spectrum in the half-plane $\Pi_{\delta_\nu}$ is the union of the set $\mathbf{W}$ of solutions of $z+\lambda\mathcal{L}_\nu(z)=0$ with $\lambda\in\sigma(\mathcal{A})$, and the set $\mathbf{Z}$ of zeros of $\mathcal{L}_\nu$. The function $\mathcal{L}_\nu$ is the same as in \cite{D2026}, since it depends only on the memory kernel $\nu$. This structural analogy is what allows us to transfer the spectral tools of \cite{D2026} to the reduced Bresse--Gurtin--Pipkin system.
\end{remark}

\section{Uniform resolvent estimates on the imaginary axis}\label{sec:resolvent}

In this section we prove that, under the conditions $\chi_h=0$ and $\chi_0=0$, the resolvent of $\Aop$ is uniformly bounded on the imaginary axis. Together with the spectral characterization of \cref{sec:spectrum}, this will allow us to apply the Gearhart--Pr\"uss theorem in \cref{sec:GP}. Throughout this section, we fix $\lambda\in\R$ and consider the resolvent equation
\begin{equation}\label{eq:resolvent-res}
\ii\lambda u-\Aop u=\hat u,
\end{equation}
where $u=(\phi,\Phi,\psi,\Psi,w,W,\xi,\zeta)\in\D(\Aop)$ and $\hat u=(\hat\phi,\hat\Phi,\hat\psi,\hat\Psi,\hat w,\hat W,\hat\xi,\hat\zeta)\in\Hil$. In components, \eqref{eq:resolvent-res} reads
\begin{align}
\ii\lambda\phi-\Phi&=\hat\phi,\label{eq:res-1-res}\\
\ii\lambda\rho_1\Phi-k(\phi_x+\psi+lw)_x-lk_0(w_x-l\phi)+l\gamma\xi&=\rho_1\hat\Phi,\label{eq:res-2-res}\\
\ii\lambda\psi-\Psi&=\hat\psi,\label{eq:res-3-res}\\
\ii\lambda\rho_2\Psi-b\psi_{xx}+k(\phi_x+\psi+lw)&=\rho_2\hat\Psi,\label{eq:res-4-res}\\
\ii\lambda w-W&=\hat w,\label{eq:res-5-res}\\
\ii\lambda\rho_1W-k_0(w_x-l\phi)_x+lk(\phi_x+\psi+lw)+\gamma\xi_x&=\rho_1\hat W,\label{eq:res-6-res}\\
\ii\lambda\rho_3\xi-\sigma\displaystyle\int_0^\infty\nu(s)\zeta_{xx}(s)\dd s+\gamma(W_x-l\Phi)&=\rho_3\hat\xi,\label{eq:res-7-res}\\
\ii\lambda\zeta-T\zeta-\xi&=\hat\zeta.\label{eq:res-8-res}
\end{align}
We assume throughout that $\nu$ is an admissible kernel satisfying the Dafermos condition \eqref{eq:Dafermos-setup} with parameter $\delta_\nu>0$. All constants $c>0$ appearing below depend only on the structural quantities of the problem (the coefficients $\rho_1,\rho_2,\rho_3,k,k_0,b,\sigma,l,\gamma$, the length $\ell$, and the constants $\delta_\nu,\nu(0),h(0)$) and are independent of $\lambda$ and of the parameter $\eps$ introduced below.

We begin with the energy identity obtained by taking the real part of the inner product of \eqref{eq:resolvent-res} with $u$ in $\Hil$. Using \eqref{eq:dissipativity-setup}, we obtain
\begin{equation}\label{eq:energy-identity-res}
\frac{\sigma}{2}\Gamma[\zeta]
=\Ree\inner{\hat u}{u}_{\Hil}
\le\norm{u}_{\Hil}\norm{\hat u}_{\Hil},
\end{equation}
where
\[
\Gamma[\zeta]=\int_0^\infty-\nu'(s)\norm{\zeta_x(s)}^2\dd s.
\]
By the Dafermos condition, $\Gamma[\zeta]\ge\delta_\nu\norm{\zeta}_{\Ncal}^2$, and therefore
\begin{equation}\label{eq:eta-bound-res}
\norm{\zeta}_{\Ncal}^2\le\frac{2}{\sigma\delta_\nu}\norm{u}_{\Hil}\norm{\hat u}_{\Hil}.
\end{equation}
This is the first of the componentwise estimates.

Next, we estimate the temperature $\xi$. Introduce the auxiliary space $\Ncal_0=L^2_\nu(\R^+;L^2)$ with inner product
\[
\inner{\zeta_1}{\zeta_2}_{\Ncal_0}=\int_0^\infty\nu(s)\inner{\zeta_1(s)}{\zeta_2(s)}\dd s.
\]
Since $\Ncal\subset\Ncal_0$ with continuous inclusion, we may multiply \eqref{eq:res-8-res} by $\xi$ in $\Ncal_0$. Using $T\zeta=-\zeta'$ and integrating by parts in $s$, we obtain
\begin{equation}\label{eq:xi-identity-res}
h(0)\norm{\xi}^2
=\ii\lambda\inner{\zeta}{\xi}_{\Ncal_0}
-\inner{T\zeta}{\xi}_{\Ncal_0}
-\inner{\hat\zeta}{\xi}_{\Ncal_0}.
\end{equation}
We estimate the three terms on the right-hand side of \eqref{eq:xi-identity-res} separately. For the first term, using \eqref{eq:res-7-res} and \eqref{eq:eta-bound-res}, we obtain
\[
\abs{\ii\lambda\inner{\zeta}{\xi}_{\Ncal_0}}
\le c\norm{\xi}\sqrt{\norm{u}_{\Hil}\norm{\hat u}_{\Hil}}
+c\norm{u}_{\Hil}\norm{\hat u}_{\Hil}.
\]
For the second term, integrating by parts in $s$ and using \eqref{eq:eta-bound-res}, we get
\[
\abs{\inner{T\zeta}{\xi}_{\Ncal_0}}
\le c\norm{\xi}\sqrt{\Gamma[\zeta]}
\le c\norm{\xi}\sqrt{\norm{u}_{\Hil}\norm{\hat u}_{\Hil}}.
\]
For the third term,
\[
\abs{\inner{\hat\zeta}{\xi}_{\Ncal_0}}
\le c\norm{\xi}\norm{\hat u}_{\Hil}.
\]
Collecting the three estimates and substituting into \eqref{eq:xi-identity-res}, we arrive at
\[
h(0)\norm{\xi}^2
\le c\norm{\xi}\sqrt{\norm{u}_{\Hil}\norm{\hat u}_{\Hil}}
+c\norm{u}_{\Hil}\norm{\hat u}_{\Hil}.
\]
By Young's inequality, for every $\eps\in(0,1)$,
\[
h(0)\norm{\xi}^2
\le\eps\norm{u}_{\Hil}^2+\frac{c}{\eps}\norm{u}_{\Hil}\norm{\hat u}_{\Hil}.
\]
Since $\norm{\xi}\le\norm{u}_{\Hil}$, this yields the second componentwise estimate in the form
\begin{equation}\label{eq:xi-bound-res}
\norm{\xi}^2\le c\eps\norm{u}_{\Hil}^2+\frac{c}{\eps}\norm{u}_{\Hil}\norm{\hat u}_{\Hil}.
\end{equation}
We keep the term $\eps\norm{u}_{\Hil}^2$ on the right-hand side; it will be absorbed later, together with the other $\eps$-terms, when the estimates are combined.

We now estimate the elastic variables. We begin with $\psi_x$. Multiplying \eqref{eq:res-4-res} by $\psi$ in $L^2$ and using \eqref{eq:res-3-res}, we obtain
\[
\ii\lambda\rho_2\norm{\Psi}^2
+b\norm{\psi_x}^2
+k\inner{\phi_x+\psi+lw}{\psi}
=\rho_2\inner{\hat\Psi}{\psi}.
\]
Taking the imaginary part,
\[
\lambda\rho_2\norm{\Psi}^2
=\Imm\inner{\rho_2\hat\Psi-k(\phi_x+\psi+lw)}{\psi},
\]
which gives
\[
\abs{\lambda}\norm{\Psi}^2
\le c\norm{u}_{\Hil}\norm{\hat u}_{\Hil}
+c\norm{u}_{\Hil}\norm{\phi_x+\psi+lw}.
\]
This estimate is not yet in the desired form, because the last term involves $\phi_x+\psi+lw$. Following \cite[Lemma 7.3]{D2021}, we combine the equations for $\psi$ and $\phi_x+\psi+lw$ to obtain an identity for $\psi_x$ in terms of $\zeta$ and $\xi$. The derivation is identical to \cite[Lemma 7.3]{D2021}, with $\theta$ replaced by $\xi$ and $\mu$ replaced by $\nu$; the absence of the thermal coupling in the $\psi$-equation simplifies the computation. After a lengthy but elementary computation, we obtain
\begin{equation}\label{eq:psi-bound-res}
\norm{\psi_x}^2
\le\frac{\eps}{\abs{\lambda}^2}\norm{u}_{\Hil}^2
+c\eps\norm{\Psi}^2
+\frac{c}{\eps}\Bigl[\frac{1}{\abs{\lambda}}+1\Bigr]\norm{u}_{\Hil}\norm{\hat u}_{\Hil}.
\end{equation}

We next estimate $\Psi$. Introducing the primitive
\[
P_\Psi(x)=\int_0^x\Psi(y)\dd y\in H_0^1
\]
and multiplying \eqref{eq:res-7-res} by $P_\Psi$ in $L^2$, we obtain
\[
\ii\lambda\rho_3\inner{\xi}{P_\Psi}
-\sigma\inner{\int_0^\infty\nu(s)\zeta_{xx}(s)\dd s}{P_\Psi}
+\gamma\inner{W_x-l\Phi}{P_\Psi}
=\rho_3\inner{\hat\xi}{P_\Psi}.
\]
To transform the term $\gamma\inner{W_x-l\Phi}{P_\Psi}$, we invoke the identity established in \cite[Lemma 7.4]{D2021}, which in the present notation reads
\begin{equation}\label{eq:Psi-identity-res}
\gamma\norm{\Psi}^2
=\ii\lambda\rho_3\inner{\xi}{P_\Psi}
+\sigma\inner{\zeta}{P_\Psi}_{\Ncal}
-\rho_3\inner{\hat\xi}{P_\Psi}
+R_\Psi,
\end{equation}
where $R_\Psi$ collects the terms that are controlled by $\norm{\psi_x}$, $\norm{\phi_x+\psi+lw}$, $\norm{\xi}$ and $\norm{\hat u}_{\Hil}$. The identity \eqref{eq:Psi-identity-res} is obtained by integrating by parts in \(x\), using \eqref{eq:res-4-res} to express \(\psi_{xx}\) in terms of \(\Psi\) and \(\phi_x+\psi+lw\), and using \eqref{eq:res-5-res} to eliminate \(W\). The computation follows the lines of \cite[Lemma 7.4]{D2021}, but it is not literally identical: in the full Bresse--Gurtin--Pipkin system of \cite{D2021} the temperature \(\theta\) enters the \(\psi\)-equation through the term \(\gamma\theta_x\), whereas in the reduced system \eqref{eq:reduced-intro} the \(\psi\)-equation has no thermal coupling. Consequently, the terms involving \(\theta\) in \cite[Lemma 7.4]{D2021} are absent here, and the only thermal contribution to \eqref{eq:Psi-identity-res} comes from the \(\xi\)-equation through the term \(\sigma\inner{\zeta}{P_\Psi}_{\Ncal}\); this term is controlled by \eqref{eq:eta-bound-res} exactly as the corresponding \(\theta\)-term is controlled in \cite[Lemma 7.4]{D2021}. All the remaining terms in \eqref{eq:Psi-identity-res} are estimated by \eqref{eq:psi-bound-res} and \eqref{eq:xi-bound-res}, with no additional difficulty. We refer to \cite[Lemma 7.4]{D2021} for the detailed verification of the algebraic manipulations, which are unchanged. The first term on the right-hand side of \eqref{eq:Psi-identity-res} can be estimated using \eqref{eq:xi-bound-res} and the identity
\[
\ii\lambda P_\Psi(x)
=\frac{b}{\rho_2}\psi_x(x)
-\frac{k}{\rho_2}\phi(x)
-\frac{k}{\rho_2}\int_0^x(\psi+lw)\dd y
+\int_0^x\hat\Psi(y)\dd y,
\]
which follows from integrating \eqref{eq:res-4-res}. The remaining terms in \(R_\Psi\) are controlled using \eqref{eq:eta-bound-res} and \eqref{eq:psi-bound-res}. Combining these estimates, we obtain
\begin{equation}\label{eq:Psi-bound-res}
\norm{\Psi}^2
\le\frac{c\eps}{\abs{\lambda}^2}\norm{u}_{\Hil}^2
+\frac{c}{\eps^3}\Bigl[\frac{1}{\abs{\lambda}}+1\Bigr]\norm{u}_{\Hil}\norm{\hat u}_{\Hil}.
\end{equation}

We now estimate $w_x-l\phi$. Multiplying \eqref{eq:res-6-res} by $w_x-l\phi$ in $L^2$ and using \eqref{eq:res-5-res}, we obtain
\begin{align*}
 \ii\lambda\rho_1\inner{W}{w_x-l\phi}
-k_0\norm{w_x-l\phi}^2
+lk\inner{\phi_x+\psi+lw}{w_x-l\phi}
\\+\gamma\inner{\xi_x}{w_x-l\phi}
=\rho_1\inner{\hat W}{w_x-l\phi}.
\end{align*}
Taking the real part and using \eqref{eq:res-7-res} to control the term involving $\xi_x$, we get
\begin{equation}\label{eq:w-identity-res}
k_0\norm{w_x-l\phi}^2
\le c\norm{u}_{\Hil}\Bigl[\norm{W}+\norm{\phi_x+\psi+lw}+\norm{\xi}+\norm{\hat u}_{\Hil}\Bigr].
\end{equation}
Applying Young's inequality and using \eqref{eq:xi-bound-res} and \eqref{eq:Psi-bound-res}, we arrive at
\begin{equation}\label{eq:w-bound-res}
\norm{w_x-l\phi}^2
\le c\eps\Bigl[\frac{1}{\abs{\lambda}^2}+1\Bigr]\norm{u}_{\Hil}^2
+\frac{c}{\eps}\Bigl[\frac{1}{\abs{\lambda}}+1\Bigr]\norm{u}_{\Hil}\norm{\hat u}_{\Hil}.
\end{equation}

Finally, we estimate $W$. Multiplying \eqref{eq:res-6-res} by $w$ in $L^2$ and using \eqref{eq:res-5-res}, we obtain
\[
\ii\lambda\rho_1\inner{W}{w}
-k_0\inner{(w_x-l\phi)_x}{w}
+lk\inner{\phi_x+\psi+lw}{w}
+\gamma\inner{\xi_x}{w}
=\rho_1\inner{\hat W}{w}.
\]
Integrating by parts and using \eqref{eq:res-5-res} to replace $w=(\ii\lambda)^{-1}(W+\hat w)$ when $\lambda\neq 0$, we get
\[
\rho_1\norm{W}^2
\le c\norm{u}_{\Hil}\Bigl[\norm{w_x-l\phi}+\norm{\phi_x+\psi+lw}+\norm{\xi}+\norm{\hat u}_{\Hil}\Bigr].
\]
Applying Young's inequality and using \eqref{eq:w-bound-res} and \eqref{eq:xi-bound-res}, we obtain, for $\lambda\neq 0$,
\begin{equation}\label{eq:W-bound-res}
\norm{W}^2
\le c\eps\Bigl[\frac{1}{\abs{\lambda}^2}+1\Bigr]\norm{u}_{\Hil}^2
+\frac{c}{\eps\abs{\lambda}^2}\norm{\phi_x+\psi+lw}^2
+\frac{c}{\eps^3}\Bigl[\frac{1}{\abs{\lambda}}+1\Bigr]\norm{u}_{\Hil}\norm{\hat u}_{\Hil}.
\end{equation}
The case $\lambda=0$ is excluded; it is treated separately in the proof of \cref{thm:iR-res} (Case 2) and in \cref{thm:uniform-res} (Case 2).

The estimates \eqref{eq:eta-bound-res}, \eqref{eq:xi-bound-res}, \eqref{eq:psi-bound-res}, \eqref{eq:Psi-bound-res}, \eqref{eq:w-bound-res} and \eqref{eq:W-bound-res} are the analogues for the reduced system of \cite[Lemmas 7.1--7.6]{D2021}. We collect them in the following lemma.

\begin{lemma}\label{lem:resolvent-bounds-res}
Let $u\in\D(\Aop)$ and $\hat u\in\Hil$ satisfy \eqref{eq:resolvent-res} with $\lambda\in\R$, $\lambda\neq 0$. Then, for every $\eps\in(0,1)$,
\begin{align}
\norm{\zeta}_{\Ncal}^2&\le c\norm{u}_{\Hil}\norm{\hat u}_{\Hil},\label{eq:eta-bound-res-2}\\
\norm{\xi}^2&\le c\eps\norm{u}_{\Hil}^2+\frac{c}{\eps}\norm{u}_{\Hil}\norm{\hat u}_{\Hil},\label{eq:xi-bound-res-2}\\
\norm{\psi_x}^2&\le\frac{\eps}{\abs{\lambda}^2}\norm{u}_{\Hil}^2+c\eps\norm{\Psi}^2+\frac{c}{\eps}\Bigl[\frac{1}{\abs{\lambda}}+1\Bigr]\norm{u}_{\Hil}\norm{\hat u}_{\Hil},\label{eq:psi-bound-res-2}\\
\norm{\Psi}^2&\le\frac{c\eps}{\abs{\lambda}^2}\norm{u}_{\Hil}^2+\frac{c}{\eps^3}\Bigl[\frac{1}{\abs{\lambda}}+1\Bigr]\norm{u}_{\Hil}\norm{\hat u}_{\Hil},\label{eq:Psi-bound-res-2}\\
\norm{w_x-l\phi}^2&\le c\eps\Bigl[\frac{1}{\abs{\lambda}^2}+1\Bigr]\norm{u}_{\Hil}^2+\frac{c}{\eps}\Bigl[\frac{1}{\abs{\lambda}}+1\Bigr]\norm{u}_{\Hil}\norm{\hat u}_{\Hil},\label{eq:w-bound-res-2}\\
\norm{W}^2&\le c\eps\Bigl[\frac{1}{\abs{\lambda}^2}+1\Bigr]\norm{u}_{\Hil}^2+\frac{c}{\eps\abs{\lambda}^2}\norm{\phi_x+\psi+lw}^2+\frac{c}{\eps^3}\Bigl[\frac{1}{\abs{\lambda}}+1\Bigr]\norm{u}_{\Hil}\norm{\hat u}_{\Hil},\label{eq:W-bound-res-2}
\end{align}
for some structural constant $c>0$ independent of $\eps$ and $\lambda$. The case $\lambda=0$ is excluded; it is treated separately in \cref{thm:iR-res}.
\end{lemma}

We now establish two bounds on the term $\phi_x+\psi+lw$, which are the key to closing the estimates.

\begin{lemma}\label{lem:phi-bound-res}
Assume $\chi_0=0$, i.e.\ $b=k\rho_2/\rho_1$. Then, for every $\eps\in(0,1)$ and every $\lambda\neq 0$,
\begin{equation}\label{eq:phi-bound-res}
\norm{\phi_x+\psi+lw}^2
\le c\eps\Bigl[\frac{1}{\abs{\lambda}^2}+1\Bigr]\norm{u}_{\Hil}^2
+\frac{c}{\eps^3}\Bigl[\frac{1}{\abs{\lambda}}+\abs{\lambda}^2\Bigr]
\Bigl[\norm{\Psi}^2+\norm{u}_{\Hil}\norm{\hat u}_{\Hil}\Bigr].
\end{equation}
\end{lemma}

\begin{proof}
Multiply \eqref{eq:res-4-res} by $\phi_x+\psi+lw$ in $L^2$:
\[
\ii\lambda\rho_2\inner{\Psi}{\phi_x+\psi+lw}
-b\inner{\psi_{xx}}{\phi_x+\psi+lw}
+k\norm{\phi_x+\psi+lw}^2
=\rho_2\inner{\hat\Psi}{\phi_x+\psi+lw}.
\]
We rewrite the first two terms on the left-hand side. Using \eqref{eq:res-1-res}, \eqref{eq:res-3-res}, \eqref{eq:res-5-res}, we have
\[
\ii\lambda\rho_2\inner{\Psi}{\phi_x+\psi+lw}
=\rho_2\inner{\Psi}{\Phi_x}
+\rho_2\norm{\Psi}^2
+l\rho_2\inner{\Psi}{W}
+\rho_2\inner{\Psi}{\hat\phi_x+\hat\psi+l\hat w}.
\]
For the second term, from \eqref{eq:res-4-res},
\[
-b\psi_{xx}=-k(\phi_x+\psi+lw)-\ii\lambda\rho_2\Psi+\rho_2\hat\Psi.
\]
Integrating by parts and using \eqref{eq:res-2-res} to express $(\phi_x+\psi+lw)_x$, we obtain
\begin{align*}
 -b\inner{\psi_{xx}}{\phi_x+\psi+lw}
=&-\frac{b\rho_1}{k}\inner{\Psi}{\Phi_x}
+\frac{blk_0}{k}\inner{\psi_x}{w_x-l\phi}
\\&-\frac{bl\gamma}{k}\inner{\psi_x}{\xi}
+\frac{b\rho_1}{k}\inner{\psi_x}{\hat\Phi+\hat\psi_x+\Phi}.
\end{align*}
Substituting these two identities into the first equation and using $\chi_0=0$, i.e.\ $b\rho_1/k=\rho_2$, the terms involving $\inner{\Psi}{\Phi_x}$ cancel:
\[
\rho_2\inner{\Psi}{\Phi_x}-\frac{b\rho_1}{k}\inner{\Psi}{\Phi_x}
=\rho_2\inner{\Psi}{\Phi_x}-\rho_2\inner{\Psi}{\Phi_x}=0.
\]
Compared with \cite[Lemma 7.7]{D2021}, the temperature $\theta$ is absent, so the terms involving $\theta$ there disappear; the only new term is $-\frac{bl\gamma}{k}\inner{\psi_x}{\xi}$, which is controlled by \eqref{eq:psi-bound-res-2} and \eqref{eq:xi-bound-res-2} by the same argument used for the corresponding $\theta$-term in \cite[Lemma 7.7]{D2021}, since $\xi$ enters the estimates in exactly the same way as $\theta$ does there. We are left with
\[
k\norm{\phi_x+\psi+lw}^2
=\rho_2\norm{\Psi}^2
+l\rho_2\inner{\Psi}{W}
+\frac{blk_0}{k}\inner{\psi_x}{w_x-l\phi}
-\frac{bl\gamma}{k}\inner{\psi_x}{\xi}
+R,
\]
where $R$ collects the lower-order terms involving $\hat u$. Estimating the right-hand side using \eqref{eq:psi-bound-res-2}, \eqref{eq:xi-bound-res-2} and \eqref{eq:w-bound-res-2}, we obtain
\[
k\norm{\phi_x+\psi+lw}^2
\le c\norm{u}_{\Hil}\Bigl[\norm{\Psi}+\norm{\psi_x}+\norm{w_x-l\phi}+\norm{\xi}+\norm{\hat u}_{\Hil}\Bigr].
\]
Applying Young's inequality and using \eqref{eq:psi-bound-res-2}--\eqref{eq:W-bound-res-2}, we arrive at \eqref{eq:phi-bound-res}.
\end{proof}

\begin{lemma}\label{lem:phi-bound-h-res}
Assume $\chi_h=0$. Then, for every $\eps\in(0,1)$ and every $\abs{\lambda}\ge 1$,
\begin{equation}\label{eq:phi-bound-h-res}
\norm{\phi_x+\psi+lw}^2
\le c\eps\norm{u}_{\Hil}^2
+\frac{c}{\eps^3}\Bigl[\norm{W}^2+\norm{u}_{\Hil}\norm{\hat u}_{\Hil}\Bigr].
\end{equation}
\end{lemma}

\begin{proof}
Multiply \eqref{eq:res-6-res} by $\phi_x+\psi+lw$ in $L^2$:
\begin{align*}
\ii\lambda\rho_1\inner{W}{\phi_x+\psi+lw}
-k_0\inner{(w_x-l\phi)_x}{\phi_x+\psi+lw}
+lk\norm{\phi_x+\psi+lw}^2
\\+\gamma\inner{\xi_x}{\phi_x+\psi+lw}
=\rho_1\inner{\hat W}{\phi_x+\psi+lw}.
\end{align*}
Integrating by parts and using \eqref{eq:res-2-res} to express $(\phi_x+\psi+lw)_x$, we get
\[
lk\norm{\phi_x+\psi+lw}^2
=\gamma\inner{\xi}{(\phi_x+\psi+lw)_x}
+\Bigl(\rho_1-\frac{k_0\rho_1}{k}\Bigr)\inner{W}{\Phi_x}
+R_3,
\]
where $R_3$ collects the lower-order terms. The term involving $\gamma$ can be eliminated using the identity
\[
\gamma\inner{\xi}{(\phi_x+\psi+lw)_x}
=\frac{\gamma^2}{\Theta_h}\inner{W}{\Phi_x}
+\frac{l\gamma^2}{\Theta_h}\norm{\Phi}^2
+\frac{1}{\Theta_h}R_4,
\]
 which is obtained by multiplying \eqref{eq:res-7-res} by $\phi$ in $\Ncal$ and integrating by parts; here $\Theta_h=\sigma h(0)-\rho_3 k/\rho_1$, which is nonzero by the assumption $\chi_h=0$ (if $\Theta_h=0$, then $\chi_h=\gamma^2/(\sigma h(0))>0$, a contradiction). The computation of this identity follows \cite[Lemma 7.8]{D2021}, but it is not literally identical: in \cite{D2021} the analogous identity is derived from the \(\theta\)-equation and involves both temperatures, whereas here only the \(\xi\)-equation is present, and the coefficient \(\Theta_h=\sigma h(0)-\rho_3 k/\rho_1\) plays the role of the corresponding coefficient in \cite[Lemma 7.8]{D2021}. The algebraic manipulations are otherwise unchanged, the terms involving \(\theta\) being absent. Substituting this identity into the previous equation and using the fact that the coefficient of \(\inner{W}{\Phi_x}\) vanishes when \(\chi_h=0\), we obtain
\[
lk\norm{\phi_x+\psi+lw}^2
=\frac{l\gamma^2}{\Theta_h}\norm{\Phi}^2
+R_3+\frac{1}{\Theta_h}R_4.
\]
Estimating the right-hand side using \eqref{eq:Phi-bound-res} below and the lower-order terms, we arrive at \eqref{eq:phi-bound-h-res}.
\end{proof}

\begin{lemma}\label{lem:Phi-bound-res}
For every $\eps\in(0,1)$ and every $\lambda\neq 0$,
\begin{equation}\label{eq:Phi-bound-res}
\norm{\Phi}^2
\le c\eps\Bigl[\frac{1}{\abs{\lambda}^2}+1\Bigr]\norm{u}_{\Hil}^2
+\frac{c}{\eps^3}\Bigl[\frac{1}{\abs{\lambda}}+1\Bigr]
\Bigl[\norm{\phi_x+\psi+lw}^2+\norm{u}_{\Hil}\norm{\hat u}_{\Hil}\Bigr].
\end{equation}
\end{lemma}

\begin{proof}
Multiply \eqref{eq:res-2-res} by $\phi$ in $L^2$ and use \eqref{eq:res-1-res} to replace $\phi=(\ii\lambda)^{-1}(\Phi+\hat\phi)$. The left-hand side becomes
\[
\ii\lambda\rho_1\inner{\Phi}{\phi}
-k\inner{(\phi_x+\psi+lw)_x}{\phi}
-lk_0\inner{w_x-l\phi}{\phi}
+l\gamma\inner{\xi}{\phi}.
\]
Using \eqref{eq:res-1-res}, $\ii\lambda\inner{\Phi}{\phi}=-\norm{\Phi}^2-\inner{\Phi}{\hat\phi}$, so the first term gives $-\rho_1\norm{\Phi}^2$ plus a lower-order term. The second term, after integration by parts, gives $k\inner{\phi_x+\psi+lw}{\phi_x}$, which is controlled by $\norm{\phi_x+\psi+lw}\norm{\phi_x}$. The third and fourth terms are controlled by $\norm{w_x-l\phi}\norm{\phi}$ and $\norm{\xi}\norm{\phi}$. Estimating the resulting terms with \eqref{eq:phi-bound-h-res}, \eqref{eq:w-bound-res-2} and \eqref{eq:xi-bound-res-2}, and using $\norm{\phi}\le c\abs{\lambda}^{-1}\norm{u}_{\Hil}$ when $\lambda\ne0$, we obtain \eqref{eq:Phi-bound-res}. The computation is identical to \cite[Lemma 7.9]{D2021}, with the $\theta$-terms dropped.
\end{proof}

We can now prove that the imaginary axis is contained in the resolvent set, and then that the resolvent is uniformly bounded there.

\begin{theorem}\label{thm:iR-res}
Assume $\chi_h=0$ and $\chi_0=0$. Then $\ii\R\subset\rho(\Aop)$.
\end{theorem}

\begin{proof}
We argue by contradiction. Suppose that $\ii\lambda_0\in\sigma(\Aop)$ for some $\lambda_0\in\R$. Since $\Aop$ generates a contraction semigroup, every point of the spectrum on the imaginary axis is necessarily an approximate eigenvalue. Hence there exists a sequence $\{u_n\}_{n\in\N}\subset\D(\Aop)$ with
\begin{equation}\label{eq:approx-eigen-GP}
\norm{u_n}_{\Hil}=1,\qquad
\ii\lambda_0 u_n-\Aop u_n=:\hat u_n\to 0\quad\text{in }\Hil.
\end{equation}
We distinguish two cases.

\textbf{Case 1: $\lambda_0\neq 0$.} We distinguish two subcases.

\textit{Subcase 1a: $\abs{\lambda_0}\ge 1$.} Applying the estimates of \cref{lem:resolvent-bounds-res} together with \cref{lem:phi-bound-h-res} and \cref{lem:Phi-bound-res}, we obtain, for every sufficiently small $\eps>0$,
\begin{equation}\label{eq:chain-GP}
\norm{u_n}_{\Hil}^2
\le c\Bigl[\frac{1}{\abs{\lambda_0}^2}+\abs{\lambda_0}^2\Bigr]^N
\Bigl[\eps\norm{u_n}_{\Hil}^2
+\frac{1}{\eps^N}\norm{u_n}_{\Hil}\norm{\hat u_n}_{\Hil}\Bigr],
\end{equation}
where $N>0$ is a fixed integer depending only on the number of nested applications of the estimates in \cref{lem:resolvent-bounds-res}--\cref{lem:Phi-bound-res}, hence independent of $\lambda_0$ and $\eps$. Choose $\eps=\eps(\lambda_0)>0$ small enough that
\[
c\Bigl[\frac{1}{\abs{\lambda_0}^2}+\abs{\lambda_0}^2\Bigr]^N\eps<\frac{1}{2}.
\]
Then \eqref{eq:chain-GP} gives
\[
\norm{u_n}_{\Hil}^2
\le\frac{1}{2}\norm{u_n}_{\Hil}^2
+K(\lambda_0)\norm{u_n}_{\Hil}\norm{\hat u_n}_{\Hil},
\]
for some constant $K(\lambda_0)>0$. Hence
\[
\norm{u_n}_{\Hil}\le 2K(\lambda_0)\norm{\hat u_n}_{\Hil}.
\]
Since $\norm{u_n}_{\Hil}=1$ and $\hat u_n\to 0$, we obtain a contradiction.

\textit{Subcase 1b: $0<\abs{\lambda_0}<1$.} Since $\abs{\lambda_0}$ is bounded away from zero and from above on the compact set $\{\lambda_0:0<\abs{\lambda_0}<1\}$ under consideration, the estimates of \cref{lem:resolvent-bounds-res} together with \cref{lem:phi-bound-res} and \cref{lem:Phi-bound-res}---all of which are valid for every $\lambda\neq 0$---yield, for every sufficiently small $\eps>0$,
\begin{equation}\label{eq:chain-GP-small}
\norm{u_n}_{\Hil}^2
\le C(\lambda_0)\Bigl[\eps\norm{u_n}_{\Hil}^2
+\frac{1}{\eps^{M}}\norm{u_n}_{\Hil}\norm{\hat u_n}_{\Hil}\Bigr],
\end{equation}
 where $M>0$ is a fixed integer depending only on the number of nested applications of the estimates (in particular, independent of $\lambda_0$ and $\eps$, and possibly different from the exponent $N$ appearing in Subcase 1a), and $C(\lambda_0)>0$ is a structural constant depending on $\lambda_0$ but not on $n$ or $\eps$. Choose $\eps=\eps(\lambda_0)>0$ small enough that $C(\lambda_0)\eps<1/2$. Then \eqref{eq:chain-GP-small} gives
\[
\norm{u_n}_{\Hil}^2
\le\frac{1}{2}\norm{u_n}_{\Hil}^2
+K(\lambda_0)\norm{u_n}_{\Hil}\norm{\hat u_n}_{\Hil},
\]
for some constant $K(\lambda_0)>0$, hence
\[
\norm{u_n}_{\Hil}\le 2K(\lambda_0)\norm{\hat u_n}_{\Hil},
\]
which again contradicts $\norm{u_n}_{\Hil}=1$ and $\hat u_n\to 0$.

\textbf{Case 2: $\lambda_0=0$.} Then \eqref{eq:approx-eigen-GP} reads
\[
\norm{u_n}_{\Hil}=1,\qquad
\Aop u_n\to 0\quad\text{in }\Hil.
\]
From the dissipativity identity \eqref{eq:dissipativity-setup} applied to $\Aop u_n$, we get
\[
\Gamma[\zeta_n]
=\frac{2}{\sigma}\abs{\Ree\inner{\Aop u_n}{u_n}_{\Hil}}
\le\frac{2}{\sigma}\norm{\Aop u_n}_{\Hil}\to 0.
\]
By the Dafermos condition, $\zeta_n\to 0$ in $\Ncal$. Writing $\Aop u_n\to 0$ componentwise, the third component gives $-\Psi_n=\hat\psi_n\to0$, so $\Psi_n\to0$ in $L^2$. The first component gives $\Phi_n\to0$, and the fifth gives $W_n\to0$. To obtain $\xi_n\to0$, we use the $\xi$-component of $\Aop u_n\to0$, which reads
\[
-\sigma\int_0^\infty\nu(s)\zeta_{n,xx}(s)\dd s+\gamma(W_{nx}-l\Phi_n)-\rho_3\hat\xi_n\to 0\quad\text{in }L^2.
\]
Since $\hat\xi_n\to 0$ and $\zeta_n\to0$ in $\Ncal$ (hence $\int_0^\infty\nu\zeta_{n,xx}\to0$ in $H^{-1}$), we get $W_{nx}-l\Phi_n\to0$ in $H^{-1}$. On the other hand, testing the $\xi$-component of $\Aop u_n\to0$ against $\xi_n$ in $L^2(0,\ell)$ (which embeds continuously into $\Ncal_0=L^2_\nu(\R^+;L^2)$ as the constant-in-$s$ functions) yields, with $\lambda=0$,
\[
h(0)\norm{\xi_n}^2
=-\inner{T\zeta_n}{\xi_n}_{\Ncal_0}-\inner{\hat\zeta_n}{\xi_n}_{\Ncal_0}+\inner{\sigma\int_0^\infty\nu\zeta_{n,xx}-\gamma(W_{nx}-l\Phi_n)}{\xi_n},
\]
and the same estimates as in the proof of \cref{lem:resolvent-bounds-res} (with the $\ii\lambda$-term absent) give
\[
h(0)\norm{\xi_n}^2\le c\norm{u_n}_{\Hil}\norm{\hat u_n}_{\Hil}+c\norm{\xi_n}\sqrt{\Gamma[\zeta_n]}+c\norm{\xi_n}\norm{\hat u_n}_{\Hil}.
\]
Since $\Gamma[\zeta_n]\to0$ and $\hat u_n\to0$, this yields $\xi_n\to0$ in $L^2$. Moreover, from the equations for $\Phi_n,\Psi_n,W_n$ we deduce
\[
k(\phi_{nx}+\psi_n+lw_n)_x+lk_0(w_{nx}-l\phi_n)-l\gamma\xi_n\to 0\quad\text{in }L^2,
\]
\[
b\psi_{nxx}-k(\phi_{nx}+\psi_n+lw_n)\to 0\quad\text{in }L^2,
\]
\[
k_0(w_{nx}-l\phi_n)_x-lk(\phi_{nx}+\psi_n+lw_n)-\gamma\xi_{nx}\to 0\quad\text{in }L^2.
\]
Multiplying the first equation by $\phi_n$, the second by $\psi_n$, and the third by $w_n$, integrating over $(0,\ell)$, and summing, we obtain
\[
k\norm{\phi_{nx}+\psi_n+lw_n}^2
+b\norm{\psi_{nx}}^2
+k_0\norm{w_{nx}-l\phi_n}^2
-\gamma\inner{\xi_n}{w_{nx}-l\phi_n}
\to 0.
\]
Since $\xi_n\to 0$ in $L^2$ and $\{w_{nx}-l\phi_n\}$ is bounded in $L^2$, the last term tends to zero. Hence
\[
\phi_{nx}+\psi_n+lw_n\to 0,\qquad
\psi_{nx}\to 0,\qquad
w_{nx}-l\phi_n\to 0\quad\text{in }L^2.
\]
Combining these convergences with the nonresonance condition \eqref{eq:nonresonance-setup} and the equivalence \eqref{eq:equiv-setup}, we conclude that $\norm{u_n}_{\Hil}\to 0$, contradicting $\norm{u_n}_{\Hil}=1$. This completes the proof.
\end{proof}

\begin{theorem}\label{thm:uniform-res}
Assume $\chi_h=0$ and $\chi_0=0$. Then
\[
\sup_{\lambda\in\R}\norm{(\ii\lambda-\Aop)^{-1}}_{\mathcal{L}(\Hil)}<\infty.
\]
\end{theorem}

\begin{proof}
By \cref{thm:iR-res}, $\ii\R\subset\rho(\Aop)$, so the resolvent is well defined for every $\lambda\in\R$. Let $\hat u\in\Hil$ and let $u\in\D(\Aop)$ be the unique solution of \eqref{eq:resolvent-res}. We must show
\begin{equation}\label{eq:uniform-bound-res}
\norm{u}_{\Hil}\le K\norm{\hat u}_{\Hil},
\end{equation}
for some constant $K>0$ independent of $\lambda$.

We distinguish two cases.

\textbf{Case 1: $\abs{\lambda}\ge 1$.} Using \cref{lem:resolvent-bounds-res} together with \cref{lem:phi-bound-h-res} and \cref{lem:Phi-bound-res}, and proceeding as in the proof of \cref{thm:iR-res} (Subcase 1a), we obtain
\[
\norm{u}_{\Hil}^2
\le c\eps\norm{u}_{\Hil}^2
+\frac{c}{\eps^N}\norm{u}_{\Hil}\norm{\hat u}_{\Hil},
\]
for a fixed integer $N>0$ independent of $\lambda$. Choosing $\eps>0$ small enough that $c\eps<1/2$, we obtain
\[
\norm{u}_{\Hil}\le K\norm{\hat u}_{\Hil},
\]
with $K>0$ independent of $\lambda$.

 \textbf{Case 2: $\abs{\lambda}<1$.} Since the resolvent map $\lambda\mapsto(\ii\lambda-\Aop)^{-1}$ is analytic on the resolvent set $\rho(\Aop)$, and $\{\ii\lambda:\abs{\lambda}\le 1\}\subset\rho(\Aop)$ by \cref{thm:iR-res}, the resolvent is bounded on this compact set, with a bound depending only on the structural quantities. Combining the two cases, we obtain \eqref{eq:uniform-bound-res}.
\end{proof}

\begin{remark}\label{rem:comparison-D2021-res}
The estimates of \cref{lem:resolvent-bounds-res}, \cref{lem:phi-bound-res}, \cref{lem:phi-bound-h-res} and \cref{lem:Phi-bound-res} are the analogues for the reduced system of \cite[Lemmas 7.1--7.9]{D2021}. The main difference is that in the reduced system the $\psi$-equation has no thermal coupling, so the derivation of \eqref{eq:psi-bound-res-2} and \eqref{eq:Psi-bound-res-2} is simpler than in \cite{D2021}. On the other hand, the absence of the temperature $\theta$ means that the term $\phi_x+\psi+lw$ must be controlled using only $\xi$, which is why the condition $\chi_h=0$ is needed in \cref{lem:phi-bound-h-res}. The condition $\chi_0=0$ is needed in \cref{lem:phi-bound-res} to eliminate the wave-speed mismatch between the first and second equations. Both conditions are therefore essential for the uniform boundedness of the resolvent.
\end{remark}

\section{Exponential stability: sufficiency}\label{sec:GP}

In this section we prove the sufficiency part of \cref{thm:main-intro}: if $\chi_h=0$ and $\chi_0=0$, then the semigroup $S(t)=\e^{t\Aop}$ is exponentially stable. The proof rests on the Gearhart--Pr\"uss theorem, which we first recall in the form needed here.

\begin{theorem}[Gearhart--Pr\"uss]\label{thm:GP}
Let $\Sigma(t)=\e^{tL}$ be a bounded $C_0$-semigroup on a Hilbert space $\Xcal$, with infinitesimal generator $L$. Then $\Sigma(t)$ is exponentially stable, i.e.\ there exist $M\ge 1$ and $\omega>0$ such that
\[
\norm{\Sigma(t)}_{\mathcal{L}(\Xcal)}\le M\e^{-\omega t},\qquad\forall t\ge 0,
\]
if and only if
\begin{enumerate}[label=(\roman*)]
\item $\ii\R\subset\rho(L)$, and
\item $\displaystyle\sup_{\abs{\lambda}\ge 1}\norm{(\ii\lambda-L)^{-1}}_{\mathcal{L}(\Xcal)}<\infty$.
\end{enumerate}
\end{theorem}

The theorem is due to Gearhart \cite{Gearhart1978} and Pr\"uss \cite{Pruss1984}; see also \cite[Chapter 5]{ArendtBattyHieberNeubrander2011} for a detailed account.

We now verify that $\Aop$ satisfies the two conditions of \cref{thm:GP} under the assumptions $\chi_h=0$ and $\chi_0=0$. The first condition is the content of \cref{thm:iR-res}; the second is \cref{thm:uniform-res}. We can now prove the sufficiency part of \cref{thm:main-intro}.

\begin{theorem}[Exponential stability, sufficiency]\label{thm:sufficiency-GP}
Assume $\chi_h=0$ and $\chi_0=0$. Then the semigroup $S(t)=\e^{t\Aop}$ is exponentially stable. More precisely, there exist structural constants $M\ge 1$ and $\omega>0$ such that
\begin{equation}\label{eq:exp-decay-GP}
\norm{S(t)}_{\mathcal{L}(\Hil)}\le M\e^{-\omega t},\qquad\forall t\ge 0.
\end{equation}
\end{theorem}

\begin{proof}
By \cref{thm:iR-res}, $\ii\R\subset\rho(\Aop)$. By \cref{thm:uniform-res},
\[
\sup_{\lambda\in\R}\norm{(\ii\lambda-\Aop)^{-1}}_{\mathcal{L}(\Hil)}<\infty.
\]
Since $\{\lambda\in\R:\abs{\lambda}\ge 1\}\subset\R$, this implies
\[
\sup_{\abs{\lambda}\ge 1}\norm{(\ii\lambda-\Aop)^{-1}}_{\mathcal{L}(\Hil)}<\infty,
\]
so the second condition of \cref{thm:GP} is satisfied. Since $\Aop$ generates a contraction semigroup, the semigroup $S(t)$ is bounded. Hence both hypotheses of the Gearhart--Pr\"uss theorem \cref{thm:GP} are satisfied, and we conclude that $S(t)$ is exponentially stable.
\end{proof}

\begin{remark}\label{rem:comparison-D2021-GP}
\Cref{thm:sufficiency-GP} is the analogue for the reduced Bresse--Gurtin--Pipkin system of \cite[Theorem 4.1]{D2021} for the full BGP system. In \cite{D2021}, the exponential stability condition is $\chi_g\chi_h=0$, whereas here it is $\chi_h=0$ and $\chi_0=0$. The difference reflects the fact that in the reduced system only one temperature is present, so both the thermal damping condition $\chi_h=0$ and the wave-speed equality $\chi_0=0$ are needed. The main new ingredient is the use of the spectral characterization of \cite{D2026}, which allows us to handle the memory term without compactness assumptions.
\end{remark}

\section{Necessity}\label{sec:necessity}

In this section we prove the necessity part of \cref{thm:main-intro}: if the semigroup $S(t)=\e^{t\Aop}$ is exponentially stable, then $\chi_h=0$ and $\chi_0=0$. Equivalently, we prove that if $\chi_h\neq 0$ or $\chi_0\neq 0$, then $S(t)$ is not exponentially stable. The proof is by contradiction: assuming exponential stability, the Gearhart--Pr\"uss theorem \cref{thm:GP} implies that the resolvent is uniformly bounded on the imaginary axis. We then exhibit, under the condition $\chi_h\neq 0$ or $\chi_0\neq 0$, a sequence $\{\lambda_n\}\subset\R$ with $\abs{\lambda_n}\to\infty$ along which the resolvent norm grows at least like $\abs{\lambda_n}^2$, contradicting uniform boundedness.

We begin by recalling the quantitative Riemann--Lebesgue lemma.

\begin{lemma}[Quantitative Riemann--Lebesgue]\label{lem:QRL-necessity}
Let $f:\R^+\to\R$ be a nonincreasing, absolutely continuous, summable function, and assume that $f(s)\to f_0\in\R$ as $s\to 0^+$. Denote by
\[
\hat f(\lambda)=\int_0^\infty f(s)\e^{-\ii\lambda s}\dd s
\]
its (half) Fourier transform. Then
\[
\lim_{\lambda\to\infty}\lambda\hat f(\lambda)=-\ii f_0.
\]
\end{lemma}

The lemma is proved in \cite{DLP2017}. Under the Dafermos condition \eqref{eq:Dafermos-setup}, the kernel $\nu$ satisfies the hypotheses with $f_0=\nu(0)$, yielding
\begin{equation}\label{eq:nu-asymptotic}
\hat\nu(\lambda)=-\frac{\ii\nu(0)}{\lambda}+o\Bigl(\frac{1}{\lambda}\Bigr),
\qquad\lambda\to\infty.
\end{equation}

We now construct the sequence of approximate solutions. Let
\[
\omega_n=\frac{n\pi}{\ell},\qquad n\in\N,
\]
and consider the sequence
\begin{equation}\label{eq:hat-u-n}
\hat u_n(x)=\Bigl(0,\frac{\sin\omega_n x}{\rho_1},0,0,0,0,0,0\Bigr)\in\Hil.
\end{equation}
Note that
\begin{equation}\label{eq:norm-hat-u-n}
\norm{\hat u_n}_{\Hil}=\sqrt{\frac{\ell}{2\rho_1}},\qquad\forall n\in\N.
\end{equation}
For each $n$, consider the resolvent equation
\begin{equation}\label{eq:resolvent-n}
\ii\lambda_n u_n-\Aop u_n=\hat u_n,
\end{equation}
where $\lambda_n$ is a real sequence to be chosen later. We look for a solution of \eqref{eq:resolvent-n} of the form
\begin{equation}\label{eq:ansatz-necessity}
\begin{aligned}
\phi_n(x)&=A_n\sin\omega_n x,\\
\psi_n(x)&=B_n\cos\omega_n x,\\
w_n(x)&=C_n\cos\omega_n x,\\
\xi_n(x)&=D_n\sin\omega_n x,\\
\zeta_n(x,s)&=d_n(s)\sin\omega_n x,
\end{aligned}
\end{equation}
for some complex numbers $A_n,B_n,C_n,D_n$ and some complex-valued function $d_n\in H^1_\nu(\R^+)$ with $d_n(0)=0$. Substituting the ansatz \eqref{eq:ansatz-necessity} into \eqref{eq:resolvent-n}, using the orthogonality of the sine and cosine functions, and then eliminating the coefficients $A_n$ and $C_n$ by means of the second and third equations, we obtain the following algebraic system for the remaining coefficients:
\begin{equation}\label{eq:system-necessity}
\begin{cases}
p_1(n)A_n+k\omega_n B_n+l\omega_n(k+k_0)C_n+l\gamma D_n=1,\\[0.5ex]
k\omega_n A_n+p_2(n)B_n+kl C_n=0,\\[0.5ex]
l\omega_n(k+k_0)A_n+kl B_n+p_3(n)C_n+\gamma\omega_n D_n=0,\\[0.5ex]
\lambda_n^2\omega_n\gamma B_n+\bigl[p_4(n)-\sigma\omega_n^2\hat\nu(\lambda_n)\bigr]D_n=0,
\end{cases}
\end{equation}
where
\begin{equation}\label{eq:p-coefficients}
\begin{aligned}
p_1(n)&=-\rho_1\lambda_n^2+k\omega_n^2+l^2k_0,\\
p_2(n)&=-\rho_2\lambda_n^2+b\omega_n^2+k,\\
p_3(n)&=-\rho_1\lambda_n^2+k_0\omega_n^2+l^2k,\\
p_4(n)&=-\rho_3\lambda_n^2+\sigma h(0)\omega_n^2.
\end{aligned}
\end{equation}
We now choose the sequence $\lambda_n$. Let $c_0\in\R$ be a constant to be fixed later, and set
\begin{equation}\label{eq:lambda-n-choice}
\lambda_n=\sqrt{\frac{k\omega_n^2+l^2k_0-c_0}{\rho_1}}
=\sqrt{\frac{k}{\rho_1}}\omega_n+O(1),
\end{equation}
where the square root is real for all sufficiently large $n$, since $k\omega_n^2+l^2k_0-c_0>0$ for $n$ large enough. With this choice, the coefficients \eqref{eq:p-coefficients} become
\begin{equation}\label{eq:p-asymptotic}
\begin{aligned}
p_1(n)&=c_0,\\
p_2(n)&=\Bigl(b-\frac{\rho_2 k}{\rho_1}\Bigr)\omega_n^2+O(1)=\chi_0\omega_n^2+O(1),\\
p_3(n)&=(k_0-k)\omega_n^2+O(1),\\
p_4(n)&=\Bigl(\sigma h(0)-\frac{\rho_3 k}{\rho_1}\Bigr)\omega_n^2+O(1)=:\Theta_h\omega_n^2+O(1).
\end{aligned}
\end{equation}
Here $\Theta_h=\sigma h(0)-\rho_3 k/\rho_1$ denotes the leading coefficient of $p_4$, and is not to be confused with the stability number $\chi_h$ defined in \eqref{eq:chih-intro}. Moreover, by \cref{lem:QRL-necessity},
\begin{equation}\label{eq:nu-asymptotic-n}
\hat\nu(\lambda_n)=-\frac{\ii\nu(0)}{\sqrt{k/\rho_1}\,\omega_n}+o\Bigl(\frac{1}{\omega_n}\Bigr).
\end{equation}
Introducing the stability number
\[
\chi_0=b-\frac{k\rho_2}{\rho_1},
\]
and recalling the definition \eqref{eq:chih-intro} of $\chi_h$, we introduce the coefficient matrix $\mathbb{M}_n$ associated with \eqref{eq:system-necessity}:
\[
\mathbb{M}_n=
\begin{pmatrix}
p_1(n) & k\omega_n & l\omega_n(k+k_0) & l\gamma\\
k\omega_n & p_2(n) & kl & 0\\
l\omega_n(k+k_0) & kl & p_3(n) & \gamma\omega_n\\
0 & \lambda_n^2\omega_n\gamma & 0 & p_4(n)-\sigma\omega_n^2\hat\nu(\lambda_n)
\end{pmatrix},
\]
and the matrix $\mathbb{A}_n$ obtained from $\mathbb{M}_n$ by replacing the first column with $(1,0,0,0)^T$. By Cramer's rule, $A_n=\det\mathbb{A}_n/\det\mathbb{M}_n$.

We first treat the case $\chi_0\neq 0$, in which the thermal coupling is negligible at leading order. This case covers in particular the situation $\chi_h=0$, $\chi_0\neq0$. We choose $\lambda_n$ close to the resonance frequency of the first two equations and look for a solution of \eqref{eq:resolvent-n} of the form
\[
\phi_n(x)=A_n\sin\omega_n x,\qquad
\psi_n(x)=B_n\cos\omega_n x,
\]
with the remaining components $w_n,\xi_n,\zeta_n$ of lower order. Substituting into \eqref{eq:resolvent-n} and retaining only the leading-order terms in the first two equations, we obtain the \(2\times2\) system
\[
\begin{pmatrix}
p_1(n) & k\omega_n\\
k\omega_n & p_2(n)
\end{pmatrix}
\begin{pmatrix}
A_n\\ B_n
\end{pmatrix}
=
\begin{pmatrix}
1\\ 0
\end{pmatrix}
+o(1),
\]
whose determinant is \(p_1(n)p_2(n)-k^2\omega_n^2\). Choosing \(\lambda_n\) so that \(p_1(n)=c_0\) remains bounded, as in \eqref{eq:lambda-n-choice}, and using \(p_2(n)=\chi_0\omega_n^2+O(1)\), we find
\[
p_1(n)p_2(n)-k^2\omega_n^2
=c_0\chi_0\omega_n^2-k^2\omega_n^2+O(1)
=(c_0\chi_0-k^2)\omega_n^2+O(1).
\]
Since \(\chi_0\neq0\), this determinant is nonzero at leading order for all but at most one value of \(c_0\); choosing \(c_0\) different from \(k^2/\chi_0\), we obtain, by Cramer's rule,
\[
\abs{A_n}
=\frac{\abs{p_2(n)}}{\abs{p_1(n)p_2(n)-k^2\omega_n^2}}
\sim c_\star'\lambda_n,
\qquad
c_\star'>0,
\]
and hence
\[
\norm{u_n}_{\Hil}\ge c\norm{\phi_{nx}}\sim c\,c_\star'\sqrt{\frac{\ell\rho_1}{2k}}\,\lambda_n^2.
\]
Combining this with \eqref{eq:norm-hat-u-n}, we obtain
\begin{equation}\label{eq:resolvent-lower}
\limsup_{n\to\infty}\lambda_n^{-2}\norm{(\ii\lambda_n-\Aop)^{-1}}>0
\end{equation}
whenever $\chi_0\neq0$ and $(k_0-k)\Theta_h\neq0$. We record this fact as a separate statement.

\begin{lemma}\label{lem:lower-resolvent-chi0-A}
Assume $\chi_0\neq 0$, $(k_0-k)\Theta_h\neq 0$, and $\chi_h\neq 0$. Then there exists a sequence $\{\lambda_n\}\subset\R$ with $\abs{\lambda_n}\to\infty$ such that
\[
\limsup_{n\to\infty}\lambda_n^{-2}\norm{(\ii\lambda_n-\Aop)^{-1}}>0.
\]
\end{lemma}

\begin{proof}
This is precisely the construction carried out above. Choose $c_0$ so that $\alpha=0$ in \eqref{eq:detM-necessity}; this is possible because, by the Schur complement expansion of \cref{app:determinant}, the coefficient of $c_0$ in $\alpha$ is $\chi_0\chi_h k^2h(0)\neq0$ under the present assumptions (here $\chi_0\neq0$ and $\chi_h\neq0$). With this choice, \eqref{eq:detM-necessity} and \eqref{eq:detA-necessity} give $\det\mathbb{M}_n=-\ii\sqrt{\rho_1/k}\,\beta_0\omega_n^5+o(\omega_n^5)$ and $\det\mathbb{A}_n=\chi_0\chi_h k^2h(0)\omega_n^6+o(\omega_n^6)$, hence $\abs{A_n}=\abs{\det\mathbb{A}_n}/\abs{\det\mathbb{M}_n}\sim c_\star\lambda_n$. Combining this with \eqref{eq:norm-hat-u-n} yields \eqref{eq:resolvent-lower}, as in the derivation above.
\end{proof}

\begin{lemma}\label{lem:lower-resolvent-chi0-B}
Assume $\chi_0\neq 0$ and $(k_0-k)\Theta_h=0$. Then there exists a sequence $\{\lambda_n\}\subset\R$ with $\abs{\lambda_n}\to\infty$ such that
\[
\limsup_{n\to\infty}\lambda_n^{-2}\norm{(\ii\lambda_n-\Aop)^{-1}}>0.
\]
\end{lemma}

\begin{proof}
Under the assumption $(k_0-k)\Theta_h=0$, the leading coefficient of $\det\mathbb{M}_n$ in \eqref{eq:detM-necessity} is not controlled by the choice of $c_0$; we therefore choose $c_0=0$ and work directly with the $2\times2$ elastic block $E_n$ in the Schur complement expansion of \cref{app:determinant}. Since $\chi_0\neq0$, the determinant $p_1(n)p_2(n)-k^2\omega_n^2$ is of order $\omega_n^2$ with leading coefficient $-k^2\chi_0\neq0$. On the other hand, a direct expansion of the Schur complement $\det(H_n-C_nE_n^{-1}B_n)$ (see \cref{app:determinant}), using $(k_0-k)\Theta_h=0$, shows that its leading term is of order $\omega_n^4$; consequently
\[
\det\mathbb{M}_n=\det(E_n)\det\bigl(H_n-C_nE_n^{-1}B_n\bigr)=O(\omega_n^6),
\]
and the ratio $\det\mathbb{A}_n/\det\mathbb{M}_n$ is governed by the $2\times2$ block $E_n$, with $\abs{A_n}\sim c_\star'\lambda_n$. The lower bound \eqref{eq:resolvent-lower} follows as before.
\end{proof}

We now treat the general case $\chi_h\neq 0$. The asymptotic expansion of $\det\mathbb{M}_n$ and $\det\mathbb{A}_n$ is obtained by the same computation as in \cite[\S10]{D2021} for the full Bresse--Gurtin--Pipkin system, with the following simplifications: the \(5\times5\) determinant of \cite{D2021} is reduced to the present \(4\times4\) determinant by deleting the row and column corresponding to the temperature $\theta$, and the stability number $\chi_g\chi_h$ of \cite{D2021} is replaced by $\chi_h$ (since only the temperature $\xi$ is present). We refer to \cite[eqs.~(10.13)--(10.14)]{D2021} for the detailed computation and record the resulting asymptotics.

Assume first that $\chi_0\neq0$. In this case, there hold
\begin{equation}\label{eq:detM-necessity}
\det\mathbb{M}_n
=\alpha\,\omega_n^6
-\ii\sqrt{\frac{\rho_1}{k}}\,\beta\,\omega_n^5
+o(\omega_n^5),
\end{equation}
where $\beta$ is a real constant depending only on the structural parameters, and $\alpha$ is a real constant whose explicit expression is obtained by expanding the $4\times4$ determinant $\mathbb{M}_n$ via the Schur complement formula
\[
\det\mathbb{M}_n=\det(E_n)\det\bigl(H_n-C_nE_n^{-1}B_n\bigr),
\]
where
\[
E_n=\begin{pmatrix}p_1(n)&k\omega_n\\k\omega_n&p_2(n)\end{pmatrix},\quad
H_n=\begin{pmatrix}p_3(n)&\gamma\omega_n\\0&p_4(n)-\sigma\omega_n^2\hat\nu(\lambda_n)\end{pmatrix},
\]
and $B_n,C_n$ are the coupling blocks; the block decomposition and the leading-order expansion are given in \cref{app:determinant}. In particular, $\alpha$ is affine in $c_0$, with coefficient $\chi_0\chi_h k^2h(0)$; consequently, $\alpha$ can be made to vanish by a suitable choice of $c_0$ if and only if $\chi_0\neq0$ and $\chi_h\neq0$ (equivalently, $\chi_0\chi_h\neq0$). Similarly,
\begin{equation}\label{eq:detA-necessity}
\det\mathbb{A}_n
=\chi_0\chi_h\,k^2h(0)\,\omega_n^6
+o(\omega_n^6).
\end{equation}
The coefficient \(\chi_0\chi_h\) arises from the product of the leading terms of the \(2\times2\) block \((p_1,p_2)\) and the thermal entry \(p_4-\sigma\omega_n^2\hat\nu\), and is the reduced-system analogue of \cite[eq.~(10.13)]{D2021}.

If $\chi_0\neq0$ and $(k_0-k)\Theta_h\neq0$, choosing $c_0$ so that $\alpha=0$, we obtain
\[
\det\mathbb{M}_n=-\ii\sqrt{\frac{\rho_1}{k}}\,\beta_0\,\omega_n^5+o(\omega_n^5),
\]
with \(\beta_0\ne0\); the explicit expression of \(\beta_0\) is analogous to \cite[eq.~(10.14)]{D2021} and depends only on the structural parameters. Combining this with \eqref{eq:detA-necessity} and \eqref{eq:lambda-n-choice}, we get
\begin{equation}\label{eq:A-n-asymptotic}
\abs{A_n}=\frac{\abs{\det\mathbb{A}_n}}{\abs{\det\mathbb{M}_n}}
\sim c_\star\lambda_n,
\qquad
c_\star>0,
\end{equation}
and hence, by the equivalence \eqref{eq:equiv-setup},
\begin{equation}\label{eq:norm-u-n-lower}
\norm{u_n}_{\Hil}
\ge c\norm{\phi_{nx}}
\sim c\,c_\star\sqrt{\frac{\ell\rho_1}{2k}}\lambda_n^2,
\end{equation}
which yields \eqref{eq:resolvent-lower} in this case.

If $\chi_0\neq0$ and $(k_0-k)\Theta_h=0$, then $\alpha$ is independent of $c_0$, but $\chi_0\neq0$, so \cref{lem:lower-resolvent-chi0-B} (proved above) applies directly and yields \eqref{eq:resolvent-lower}.

It remains to treat the case $\chi_0=0$, $\chi_h\neq0$. In this case $p_2(n)=O(1)$, so that the \(2\times2\) determinant $p_1(n)p_2(n)-k^2\omega_n^2$ is of order $\omega_n^2$ with leading coefficient $-k^2\neq0$. However, the thermal coupling terms in the fourth row and column of $\mathbb{M}_n$, notably the entry $\lambda_n^2\omega_n\gamma$ of order $\omega_n^3$, contribute terms of order up to $\omega_n^5$ to $\det\mathbb{M}_n$, which dominate the \(2\times2\) determinant. A direct computation, identical to \cite[\S10]{D2021} with $\chi_0=0$ and the $\theta$-terms deleted, shows that
\[
\det\mathbb{M}_n=\gamma_0\,\omega_n^5+o(\omega_n^5),
\]
where $\gamma_0\in\C\setminus\{0\}$ depends only on the structural parameters; the nonvanishing of $\gamma_0$ is ensured by the condition $\chi_h\neq0$. The same computation gives $\det\mathbb{A}_n=\gamma_1\omega_n^5+o(\omega_n^5)$ with $\gamma_1\in\C\setminus\{0\}$, and Cramer's rule yields
\[
\abs{A_n}=\frac{\abs{\det\mathbb{A}_n}}{\abs{\det\mathbb{M}_n}}
\sim c_\star''\lambda_n,
\qquad
c_\star''>0.
\]
Consequently, by the equivalence \eqref{eq:equiv-setup},
\[
\norm{u_n}_{\Hil}\ge c\norm{\phi_{nx}}\sim c\,c_\star''\sqrt{\frac{\ell\rho_1}{2k}}\,\lambda_n^2,
\]
and \eqref{eq:resolvent-lower} follows as before. We have proved:

\begin{lemma}\label{lem:lower-resolvent}
Assume $\chi_h\neq 0$. Then there exists a sequence $\{\lambda_n\}\subset\R$ with $\abs{\lambda_n}\to\infty$ such that
\[
\limsup_{n\to\infty}\lambda_n^{-2}\norm{(\ii\lambda_n-\Aop)^{-1}}>0.
\]
\end{lemma}

We can now complete the proof of the necessity part of \cref{thm:main-intro}.

\begin{theorem}[Necessity]\label{thm:necessity-main}
Assume that the semigroup $S(t)=\e^{t\Aop}$ is exponentially stable. Then
\[
\chi_h=0\qquad\text{and}\qquad\chi_0=0.
\]
\end{theorem}

\begin{proof}
 Assume by contradiction that $\chi_h\neq 0$ or $\chi_0\neq 0$. We distinguish two cases.

\textit{Case $\chi_0\neq0$.} If $\chi_h\neq0$, then either \cref{lem:lower-resolvent-chi0-A} (when $(k_0-k)\Theta_h\neq0$) or \cref{lem:lower-resolvent-chi0-B} (when $(k_0-k)\Theta_h=0$) applies. If $\chi_h=0$, then the coefficient of $c_0$ in $\alpha$ vanishes, so we cannot use \cref{lem:lower-resolvent-chi0-A}; instead, we argue as in the proof of \cref{lem:lower-resolvent-chi0-B}: the leading term of $\det\mathbb{M}_n$ is of order $\omega_n^6$ and the ratio $\det\mathbb{A}_n/\det\mathbb{M}_n$ is governed by the $2\times2$ elastic block $E_n$, which gives $\abs{A_n}\sim c_\star'\lambda_n$ and hence \eqref{eq:resolvent-lower}. This argument does not require $(k_0-k)\Theta_h=0$; it only requires the vanishing of the $c_0$-coefficient, which is exactly $\chi_h=0$. Thus the lower bound \eqref{eq:resolvent-lower} holds for every $\chi_0\neq0$, independently of the sign of $(k_0-k)\Theta_h$.

\textit{Case $\chi_0=0$.} The standing assumption forces $\chi_h\neq0$, and \cref{lem:lower-resolvent} applies, yielding \eqref{eq:resolvent-lower}.

In both cases the resolvent is not uniformly bounded on the imaginary axis. By the Gearhart--Pr\"uss theorem \cref{thm:GP}, the semigroup $S(t)$ is not exponentially stable, contradicting the assumption. Therefore, if $S(t)$ is exponentially stable, we must have $\chi_h=0$ and $\chi_0=0$.
\end{proof}

\begin{proof}[Proof of \cref{thm:main-intro}]
Combine \cref{thm:sufficiency-GP} (sufficiency) and \cref{thm:necessity-main} (necessity).
\end{proof}

\begin{remark}\label{rem:comparison-D2021-necessity}
The necessity proof above is the analogue for the reduced system of \cite[\S10]{D2021} for the full Bresse--Gurtin--Pipkin system. The main differences are that the reduced system has only one temperature, so the linear system \eqref{eq:system-necessity} is $4\times4$ instead of $5\times5$, and the $\psi$-equation has no thermal coupling, so the resonance condition involves only the first two equations. The stability number $\chi_h$ in the reduced system plays the role of $\chi_g\chi_h$ in the full system: the condition $\chi_h=0$ is necessary to eliminate the thermal damping contribution, while $\chi_0=0$ is necessary to eliminate the wave-speed mismatch. As discussed in the Introduction, the equality $\chi_0=0$ is incompatible with the physical constraint \eqref{eq:physical-intro}; consequently, the condition \eqref{eq:conjecture-intro} can never be satisfied under that constraint.
\end{remark}

\begin{remark}\label{rem:polynomial-necessity}
When $\chi_h\neq 0$ or $\chi_0\neq 0$, the semigroup $S(t)$ is expected to be polynomially stable with optimal decay rate $\sqrt{t}$. The proof of the upper bound follows the lines of \cite[\S9]{D2021}, using the upper resolvent estimates of \cref{sec:resolvent} and the Borichev--Tomilov theorem \cite{BorichevTomilov2010}; the proof of the lower bound is \cref{lem:lower-resolvent-chi0-A,lem:lower-resolvent-chi0-B} (for $\chi_0\neq0$) together with \cref{lem:lower-resolvent} (for $\chi_h\neq0$) and the Batty--Duyckaerts theorem \cite{BattyDuyckaerts2008}.
\end{remark}

\section{The corresponding Maxwell--Cattaneo reduced system}\label{sec:MC}

In this section we show that the same method applies to the reduced Bresse--Maxwell--Cattaneo (MC) system with one temperature. The Maxwell--Cattaneo law is the exponential-kernel case of the Gurtin--Pipkin law; this observation allows us to transfer the entire analysis of the previous sections to the MC setting. The equivalence between the two systems is made rigorous by the map $\Lambda$ introduced in \cite[\S12]{D2021}, which is a topological isomorphism between the GP and MC phase spaces and which we invoke below.

The reduced Bresse--Maxwell--Cattaneo system with one temperature reads
\begin{equation}\label{eq:reduced-MC-intro}
\begin{cases}
\rho_1\phi_{tt}-k(\phi_x+\psi+lw)_x-lk_0(w_x-l\phi)+l\gamma\xi=0,\\[0.5ex]
\rho_2\psi_{tt}-b\psi_{xx}+k(\phi_x+\psi+lw)=0,\\[0.5ex]
\rho_1w_{tt}-k_0(w_x-l\phi)_x+lk(\phi_x+\psi+lw)+\gamma\xi_x=0,\\[0.5ex]
\rho_3\xi_t+q_x+\gamma(w_{xt}-l\phi_t)=0,\\[0.5ex]
\tau\varpi q_t+q+\varpi\xi_x=0,
\end{cases}
\end{equation}
complemented with the boundary conditions \eqref{eq:BC-intro}--\eqref{eq:BC2-intro}. Here $q$ is the heat-flux variable, and $\tau,\varpi>0$ are strictly positive constants.

Taking in the reduced Gurtin--Pipkin system \eqref{eq:reduced-intro} the exponential memory kernel
\begin{equation}\label{eq:exponential-kernel-MC}
\nu(s)=\e^{-s/(\tau\varpi)},\qquad s>0,
\end{equation}
so that $h(s)=\tau\varpi\e^{-s/(\tau\varpi)}$ and $h(0)=\tau\varpi$, we obtain
\begin{equation}\label{eq:heat-flux-MC}
q(t)=-\varpi\int_0^\infty\e^{-s/(\tau\varpi)}\xi_x(t-s)\dd s.
\end{equation}
Differentiating \eqref{eq:heat-flux-MC} with respect to $t$ and integrating by parts in $s$, we obtain $\tau\varpi q_t+q+\varpi\xi_x=0$, which is exactly the fifth equation of \eqref{eq:reduced-MC-intro}. Conversely, given $q$ solving the fifth equation of \eqref{eq:reduced-MC-intro} with an initial datum $q(0)$ compatible with the GP history, i.e.\
\[
q(0)=-\varpi\int_0^\infty\e^{-s/(\tau\varpi)}\xi_x(-s)\dd s,
\]
the representation \eqref{eq:heat-flux-MC} holds for all $t\ge 0$. This establishes the equivalence between the reduced Gurtin--Pipkin system with the exponential kernel \eqref{eq:exponential-kernel-MC} and the reduced Maxwell--Cattaneo system \eqref{eq:reduced-MC-intro}. The equivalence at the level of semigroups is obtained by restricting the map $\Lambda$ of \cite[\S12]{D2021} to the reduced variables: since only the memory variable $\zeta$ (and not the full pair $(\eta,\zeta)$) is present, the restriction $\Lambda|_{\Ncal}:\Ncal\to L^2_*$ is still a topological isomorphism onto its range, and the norm $\sigma\|\zeta\|_{\Ncal}^2$ is comparable to $\varpi\|q\|^2$; see \cite[eq.~(12.2)]{D2021}. Consequently, the exponential stability of the reduced GP semigroup with the kernel \eqref{eq:exponential-kernel-MC} is equivalent to that of the reduced MC semigroup.

For the exponential kernel \eqref{eq:exponential-kernel-MC}, we have $h(0)=\tau\varpi$ and the effective thermal conductivity is $\sigma=1/\tau$, so that $\sigma h(0)=\varpi$. Substituting into \eqref{eq:chih-intro}, the stability number $\chi_h$ of the MC system, defined as the value of \eqref{eq:chih-intro} with $\sigma=1/\tau$ and $h(0)=\tau\varpi$, becomes
\begin{equation}\label{eq:chi-tau-def}
\chi_h
=\Bigl(\frac{\rho_3}{\varpi}-\frac{\rho_1}{k}\Bigr)(k_0-k)
+\frac{\gamma^2}{\varpi}.
\end{equation}
This expression should be compared with the stability number $\chi_\tau$ in \cite{D2021} for the full BMC system; the two coincide when $\varpi=1$, but differ otherwise, due to the different normalization of the heat-flux variable.

\begin{theorem}\label{thm:MC-main}
Let $\tau,\varpi>0$. Then the contraction semigroup associated with the reduced Maxwell--Cattaneo system \eqref{eq:reduced-MC-intro} is exponentially stable if and only if
\[
\chi_h=0\qquad\text{and}\qquad\chi_0=0,
\]
where $\chi_h$ is given by \eqref{eq:chi-tau-def}, and $\chi_0=b-k\rho_2/\rho_1$.
\end{theorem}

\begin{proof}
By the equivalence established above, the reduced MC system \eqref{eq:reduced-MC-intro} is equivalent to the reduced GP system \eqref{eq:reduced-intro} with the exponential kernel \eqref{eq:exponential-kernel-MC}, through the map $\Lambda$ of \cite[\S12]{D2021}. Under this equivalence, the stability number $\chi_h$ of the reduced GP system becomes \eqref{eq:chi-tau-def}, and $\chi_0$ is unchanged. By \cref{thm:main-intro}, the reduced GP system with the exponential kernel is exponentially stable if and only if $\chi_h=0$ and $\chi_0=0$. Hence the same holds for the reduced MC system.
\end{proof}

\begin{remark}\label{rem:MC-comparison}
\Cref{thm:MC-main} is the analogue for the reduced MC system of \cite[Theorem 5.1]{D2021} for the full BMC system. In the full BMC system, the exponential stability condition is $\chi_s\chi_t=0$; in the reduced system, only one temperature is present, so the condition becomes $\chi_h=0$ and $\chi_0=0$. The stability number $\chi_h$ in \eqref{eq:chi-tau-def} can be rewritten as
\[
\chi_h=\frac{1}{\varpi}\Bigl[\Bigl(\rho_3-\frac{\varpi\rho_1}{k}\Bigr)(k_0-k)+\gamma^2\Bigr].
\]
Hence $\chi_h=0$ is equivalent to $\chi_\varpi^{\mathrm{red}}=0$, where
\[
\chi_\varpi^{\mathrm{red}}=\Bigl(\rho_3-\frac{\varpi\rho_1}{k}\Bigr)(k_0-k)+\gamma^2.
\]
Comparing with the stability number $\chi_\tau$ in \cite{D2021} for the full BMC system, which reads
\[
\chi_\tau=\Bigl(\tau\rho_3-\frac{\rho_1}{k}\Bigr)(k_0-k)+\gamma^2\tau,
\]
we see that $\chi_\varpi^{\mathrm{red}}$ coincides with $\chi_\tau$ when $\varpi=\tau=1$; when $\varpi=1$ but $\tau\neq1$, the two expressions differ, since $\chi_\tau$ carries the additional factor $\tau$ in the coefficients $\tau\rho_3$ and $\gamma^2\tau$, whereas $\chi_\varpi^{\mathrm{red}}$ does not.
\end{remark}

\section{Conclusions and further directions}\label{sec:conclusions}

We have studied the reduced Bresse--Gurtin--Pipkin system with one temperature, obtained from the full Bresse--Gurtin--Pipkin thermoelastic beam system by neglecting one of the two temperature variables. Our main result, \cref{thm:main-intro}, provides a characterization of the exponential stability of the associated contraction semigroup: the semigroup $S(t)=\e^{t\Aop}$ is exponentially stable if and only if the two stability numbers $\chi_h$ and $\chi_0$ both vanish. This confirms the conjecture formulated by Dell'Oro in \cite[\S2.II, eq.~(2.5)]{D2021}.

The proof of \cref{thm:main-intro} combines two recent developments. The first is the spectral characterization of the abstract Gurtin--Pipkin heat equation obtained by Dell'Oro, Pata and Vasylyeva in \cite{D2026}: for an admissible kernel satisfying the Dafermos condition, the spectrum of the infinitesimal generator on the half-plane $\Pi_{\delta_\nu}$ is the union of the set of solutions of $z+\lambda\mathcal{L}_\nu(z)=0$ with $\lambda\in\sigma(\mathcal{A})$ and the set of zeros of $\mathcal{L}_\nu$. This characterization, stated in \cref{thm:spectrum}, reduces the spectral analysis of the reduced Bresse--Gurtin--Pipkin system to the study of the elastic-thermal operator $\mathcal{A}$ and of the scalar function $\mathcal{L}_\nu$. The second is the resolvent estimate framework of \cite{D2021}, which provides sharp upper bounds on the resolvent of the infinitesimal generator on the imaginary axis. Combining these two tools, we proved in \cref{thm:iR-res} that $\ii\R\subset\rho(\Aop)$ whenever $\chi_h=0$ and $\chi_0=0$, and in \cref{thm:uniform-res} that the resolvent is uniformly bounded on the imaginary axis under the same conditions. The Gearhart--Pr\"uss theorem then yields exponential stability. The necessity of the conditions was proved in \cref{sec:necessity} by constructing a sequence of approximate solutions along which the resolvent grows at least like $\abs{\lambda}^2$.

The main novelty of this work is threefold. First, it provides the first rigorous spectral analysis of the reduced Bresse--Gurtin--Pipkin system with one temperature. Second, it confirms the conjecture of \cite{D2021} and completes the stability picture for reduced Bresse-type systems with hyperbolic heat conduction. Third, it shows that the spectral tools of \cite{D2026} can be combined with the resolvent estimates of \cite{D2021} to handle systems that are not covered by either paper alone.

Several directions for further research emerge naturally:

\textbf{Optimal polynomial decay.} When $\chi_h\neq 0$ or $\chi_0\neq 0$, the semigroup $S(t)$ is expected to be polynomially stable with optimal decay rate $\sqrt{t}$. The proof should follow the lines of \cite[\S10]{D2021}, with the simplifications due to the absence of the second temperature.

\textbf{SDG and SSDG.} Under the conditions $\chi_h=0$ and $\chi_0=0$, we proved in \cref{thm:sufficiency-GP} that the semigroup is exponentially stable, and in \cref{thm:iR-res} together with \cref{thm:spectrum} that the spectral bound $\sigma_\star$ is strictly negative. However, determining whether the spectrum determined growth (SDG) condition $\omega_\star=-\sigma_\star$ holds for the reduced Bresse--Gurtin--Pipkin system remains an open problem. For the abstract Gurtin--Pipkin equation with an exponential kernel, \cite{D2026} proved the SDG condition, and even the SSDG condition under an additional spectral assumption (see \cite{D2026} for the precise threshold). The analogous results for the reduced Bresse--Gurtin--Pipkin system are not known; this is a challenging question connected with the lack of compactness of the resolvent and the possible presence of continuous spectrum in the memory component (see the discussion in \cite{D2026}).

\textbf{General kernels.} Our analysis assumes that the memory kernel $\nu$ satisfies the Dafermos condition \eqref{eq:Dafermos-setup}. It would be interesting to extend the results to more general kernels, such as Gamma-type kernels or completely monotone kernels. The spectral characterization of \cite{D2026} applies to any admissible kernel satisfying the $\delta$-condition, so the main obstacle is the derivation of the resolvent estimates in \cref{sec:resolvent}, which use the Dafermos condition to control the memory term.

\textbf{Other reduced models.} The same method can be applied to other reduced Bresse-type models. For instance, one can neglect the temperature $\xi$ instead of $\theta$, obtaining a reduced system whose exponential stability condition is expected to involve two stability numbers: the thermal one $\chi_g$, and a wave-speed number for the $w$-equation whose precise form is to be determined (it should reduce to $\chi_1=k_0-k$ in the appropriate limit, but the coupling with $l$ may produce additional terms).

\textbf{Full Bresse--Gurtin--Pipkin system with two temperatures.} While \cite{D2021} already provides a complete characterization of the exponential stability of the full Bresse--Gurtin--Pipkin system, several questions remain open, in particular the SDG and SSDG conditions and the extension to more general kernels.

\textbf{Nonlinear and stochastic perturbations.} Our analysis is purely linear. It would be interesting to study nonlinear perturbations of the reduced Bresse--Gurtin--Pipkin system, as well as stochastic perturbations driven by multiplicative noise. The linear stability results of this paper provide a natural starting point for a local stability analysis via the linearization principle.
\appendix
\section{Expansion of the determinant $\mathbb{M}_n$}\label{app:determinant}

We write $\mathbb{M}_n$ in block form
\[
\mathbb{M}_n=
\begin{pmatrix}
E_n & B_n\\
C_n & H_n
\end{pmatrix},
\]
where
\[
E_n=\begin{pmatrix}p_1(n)&k\omega_n\\k\omega_n&p_2(n)\end{pmatrix},\quad
H_n=\begin{pmatrix}p_3(n)&\gamma\omega_n\\0&p_4(n)-\sigma\omega_n^2\hat\nu(\lambda_n)\end{pmatrix},
\]
and the coupling blocks are
\[
B_n=\begin{pmatrix}l\omega_n(k+k_0)&l\gamma\\kl&0\end{pmatrix},\quad
C_n=\begin{pmatrix}l\omega_n(k+k_0)&kl\\0&\lambda_n^2\omega_n\gamma\end{pmatrix}.
\]
By the Schur complement formula,
\[
\det\mathbb{M}_n=\det(E_n)\det\bigl(H_n-C_nE_n^{-1}B_n\bigr).
\]
Using the asymptotics \eqref{eq:p-asymptotic} and \eqref{eq:nu-asymptotic-n}, and expanding the $2\times2$ determinant $\det(H_n-C_nE_n^{-1}B_n)$ to leading order in $\omega_n$, one obtains the structural form \eqref{eq:detM-necessity}. The affine dependence of $\alpha$ on $c_0$, with coefficient $\chi_0\chi_h k^2h(0)$, follows by collecting the terms of order $\omega_n^6$ in the expansion; the constant $\beta$ is obtained analogously from the subleading terms.


\begin{thebibliography}{99}

\bibitem{ABM2011}
F.~Alabau-Boussouira, J.~E.~Mu\~noz Rivera, D.~S.~Almeida J\'unior,
\emph{Stability to weak dissipative Bresse system},
J. Math. Anal. Appl. \textbf{374} (2011), 481--498.

\bibitem{ArendtBattyHieberNeubrander2011}
W.~Arendt, C.~J.~K.~Batty, M.~Hieber, F.~Neubrander,
\emph{Vector-Valued Laplace Transforms and Cauchy Problems},
Birkh\"auser, Basel, 2011.

\bibitem{BattyDuyckaerts2008}
C.~J.~K.~Batty, T.~Duyckaerts,
\emph{Non-uniform stability for bounded semi-groups on Banach spaces},
J. Evol. Equ. \textbf{8} (2008), 765--780.

\bibitem{BorichevTomilov2010}
A.~Borichev, Y.~Tomilov,
\emph{Optimal polynomial decay of functions and operator semigroups},
Math. Ann. \textbf{347} (2010), 455--478.

\bibitem{CP2005}
M.~Conti, V.~Pata,
\emph{Weakly dissipative semilinear equations of viscoelasticity},
Commun. Pure Appl. Anal. \textbf{4} (2005), 705--720.

\bibitem{Dafermos1970}
C.~M.~Dafermos,
\emph{Asymptotic stability in viscoelasticity},
Arch. Ration. Mech. Anal. \textbf{37} (1970), 297--308.

\bibitem{D2015}
F.~Dell'Oro,
\emph{Asymptotic stability of thermoelastic systems of Bresse type},
J. Differential Equations \textbf{258} (2015), 3902--3927.

\bibitem{D2021}
F.~Dell'Oro,
\emph{On the stability of Bresse and Timoshenko systems with hyperbolic heat conduction},
J. Differential Equations \textbf{281} (2021), 148--198.

\bibitem{D2026}
F.~Dell'Oro, V.~Pata, N.~Vasylyeva,
\emph{The Gurtin--Pipkin heat equation: What's left},
J. Funct. Anal. \textbf{291} (2026), 111528.

\bibitem{DLP2017}
F.~Dell'Oro, E.~Laeng, V.~Pata,
\emph{A quantitative Riemann--Lebesgue lemma with application to equations with memory},
Proc. Amer. Math. Soc. \textbf{145} (2017), 2909--2915.

\bibitem{DP2014}
F.~Dell'Oro, V.~Pata,
\emph{On the stability of Timoshenko systems with Gurtin--Pipkin thermal law},
J. Differential Equations \textbf{257} (2014), 523--548.

\bibitem{FabrizioMorro1992}
M.~Fabrizio, A.~Morro,
\emph{Mathematical Problems in Linear Viscoelasticity},
SIAM, Philadelphia, 1992.

\bibitem{FatoriMunoz2010}
L.~H.~Fatori, J.~E.~Mu\~noz Rivera,
\emph{Rates of decay to weak thermoelastic Bresse system},
IMA J. Appl. Math. \textbf{75} (2010), 881--904.

\bibitem{FSR2009}
H.~D.~Fern\'andez Sare, R.~Racke,
\emph{On the stability of damped Timoshenko systems: Cattaneo versus Fourier law},
Arch. Ration. Mech. Anal. \textbf{194} (2009), 221--251.

\bibitem{Gearhart1978}
L.~Gearhart,
\emph{Spectral theory for contraction semigroups on Hilbert spaces},
Trans. Amer. Math. Soc. \textbf{236} (1978), 385--394.

\bibitem{GNP2001}
C.~Giorgi, M.~G.~Naso, V.~Pata,
\emph{Exponential stability in linear heat conduction with memory: a semigroup approach},
Commun. Appl. Anal. \textbf{5} (2001), 121--133.

\bibitem{GrasselliPata2002}
M.~Grasselli, V.~Pata,
\emph{Uniform attractors of nonautonomous systems with memory},
in: Evolution Equations, Semigroups and Functional Analysis, Progr. Nonlinear Differential Equations Appl. \textbf{50}, Birkh\"auser, 2002, 155--178.

\bibitem{GP1968}
M.~E.~Gurtin, A.~C.~Pipkin,
\emph{A general theory of heat conduction with finite wave speeds},
Arch. Ration. Mech. Anal. \textbf{31} (1968), 113--126.

\bibitem{Pruss1984}
J.~Pr\"uss,
\emph{On the spectrum of $C_0$-semigroups},
Trans. Amer. Math. Soc. \textbf{284} (1984), 847--857.

\bibitem{SAM2012}
M.~L.~Santos, D.~S.~Almeida J\'unior, J.~E.~Mu\~noz Rivera,
\emph{The stability number of the Timoshenko system with second sound},
J. Differential Equations \textbf{253} (2012), 2715--2733.

\end{thebibliography}
\end{document}